\documentclass[12pt]{article}
\usepackage{amsmath,amsthm,amsfonts,amssymb,amscd, amsxtra}
\theoremstyle{plain}
\newtheorem{theorem}{Theorem}[section]
\newtheorem{lemma}{Lemma}[section]
\newtheorem{definition}{Definition}[section]
\newtheorem{corollary}{Corollary}[section]
\newtheorem{proposition}{Proposition}[section]

\newtheorem{remark}{Remark}[section]
\newtheorem{claim}{Claim}[section]

\newcommand{\beq}{\begin{equation}}
\newcommand{\eeq}{\end{equation}}
\newcommand{\beqa}{\begin{eqnarray}}
\newcommand{\eeqa}{\end{eqnarray}}
\newcommand{\beqas}{\begin{eqnarray*}}
\newcommand{\eeqas}{\end{eqnarray*}}

\def\argmin{\operatorname{argmin}}

\def\min{\operatorname{min}}
\def\max{\operatorname{max}}
\DeclareMathOperator{\tr}{tr}
\DeclareMathOperator{\grad}{grad}
\DeclareMathOperator{\conv}{conv}

\DeclareMathOperator{\inter}{int}
\DeclareMathOperator{\diag}{diag}
\begin{document}

\title{Proximal Point Method for a Special Class of Nonconvex \\Functions
on Hadamard Manifolds}

\author{
G. C. Bento\thanks{IME, Universidade Federal de Goi\'as,
Goi\^ania, GO 74001-970, BR (Email: {\tt glaydston@mat.ufg.br})}
\and
O. P. Ferreira
\thanks{IME, Universidade Federal de Goi\'as,
Goi\^ania, GO 74001-970, BR (Email: {\tt orizon@mat.ufg.br}).
The author was supported in part by CNPq Grant 302618/2005-8,  PRONEX--Optimization(FAPERJ/CNPq) and FUNAPE/UFG.}\and
P. R. Oliveira \thanks{COPPE-Sistemas, Universidade Federal do Rio de Janeiro,
Rio de Janeiro, RJ 21945-970, BR (Email: {\tt poliveir@cos.ufrj.br}).
This author was supported in part by CNPq.}
}
\date{April 13, 2010}

\maketitle
\begin{abstract}
In this paper we present the proximal point method for a special class of nonconvex function on a Hadamard manifold. The well definedness of the sequence generated by the proximal point method is guaranteed. Moreover, it is proved that each accumulation point of this sequence satisfies the necessary optimality conditions and, under additional assumptions, its convergence for a minimizer is obtained.
\end{abstract}
\section{Introduction}
The extension of the concepts and techniques of the Mathematical Programming from Euclidean space to Riemannian manifolds occurs naturally and has been frequently done in recent years, with a theoretical purpose as well as to obtain effective algorithms; see \cite{Absil2007}, \cite{Teboulle2004}, \cite{BP07}, \cite{FO06}, \cite{FO98}, \cite{FO2000}, \cite{Zhu2007}, \cite{Chong Li2009}, \cite{Chong Li2009-2}, \cite{PP08}, \cite{PO2009}, \cite{RAP97}, \cite{S94} and \cite{U94}. In particular, we observe that, these extensions allow the solving of some nonconvex constrained problems in Euclidean space. More precisely, nonconvex problems in the classic sense may become convex with the introduction of an adequate Riemannian metric on the manifold (see, for example \cite{XFLN2006}). The proximal point algorithm, introduced by Martinet \cite{M70} and Rockafellar \cite{R76}, has been extended to different contexts, see \cite{FO2000}, \cite{Chong Li2009}, \cite{PO2009} and their references. In \cite{FO2000} the proximal point method has been generalized in order to solve convex optimization problems of the form
\begin{eqnarray}\label{po:conv}
\begin{array}{clc}
  (P) & \min  f(p) \\
   & \textnormal{s.t.}\,\,\, p\in M,\\
\end{array}
\end{eqnarray}
where $M$ is a Hadamard manifold and $f:M\to\mathbb{R}$ is a convex function (in Riemannian sense). The method was described as follows:
\begin{equation}\label{E:1.2}
p^{k+1}:=\argmin_{p\in M} \left\{f(p)+\frac{\lambda_k}{2}d^2(p,{p^k})\right\},
\end{equation}
with $p^{\circ}\in M$ an arbitrary point, $d$ the intrinsic Riemannian distance (to be defined later on) and $\{\lambda_k\}$ a sequence of positive numbers. The authors also showed that this extension is natural. In \cite{Chong Li2009} the important notion of maximal monotonicity from a multivalued operator defined on a Banach space to multivalued vector field defined on a Hadamard manifold has been extended. Beside the authors present a general proximal point method to finding singularity of a multivalued vector field. In particular, as an application of the convergence result obtained for the proposed algorithm,  constrained optimization problems have been solved. With regards to \cite{PO2009} the authors generalized the proximal point method with Bregman distance for solving quasiconvex and convex optimization problems on Hadamard manifolds.  Spingarn has in \cite{Spingarn1982}, in particular, developed the proximal point method for the minimization of a certain class of nondifferentiable and noncovex functions, namely, lower-$C^2$ functions defined on the Euclidean space, see also \cite{Sagastizabal2005}. Kaplan and Tichatschke in \cite{Kaplan1998} also applied the proximal point method for the minimization of a similar class of the ones studied in \cite{Sagastizabal2005} and \cite{Spingarn1982}, namely, the maximum of continuously differentiable functions.

Our goal is to study the same class objective functions studied by Kaplan and Tichatschke, \cite{Kaplan1998}, in the Riemannian context, applying the proximal point method \eqref{E:1.2} in order to solve the problem \eqref{po:conv} with the objective function in that class. For this purpose, it is necessary to study the generalized directional derivative and subdifferential in the Riemannian manifolds context. Several works have studied such concepts and presented many useful results in the Riemannian nonsmooth optimization context, see for example \cite{Azagra2005}, \cite{Zhu2007}, \cite{Montreanu-Pavel1982} and \cite{Thãmelt1992}.

The organization of our paper is as follows. In Section \ref{sec2} we define the notations and list some results of Riemannian geometry to be used throughout this paper. In Section \ref{sec3}, we recall some facts of the convex analysis on Hadamard manifolds. In Section \ref{sec:dd} we present some properties of the directional derivative of a convex function defined on a Hadamard manifold, including a characterization of the directional derivative and of the subdifferential of the maximum of a certain class of convex functions. Here we also present the definition for the generalized directional derivative of locally Lipschitz functions (not necessarily convex) and an important property of the  subdifferential of the maximum of differentiable continuously functions. In Section \ref{sec5} we present an application of the proximal point method (\ref{E:1.2}) to solve the problem \eqref{po:conv} in the case where the objective function is a real-valued function on a Hadamard manifold $M$ (not necessarily convex) given by the maximum of a certain class of functions. The main results are the proof of well definition of the sequence generated by (\ref{E:1.2}), the proof that each accumulation point of this sequence is a stationary point of the objective function and, under some additional assumptions, the proof of convergence of that sequence to a solution of the problem \eqref{po:conv}. Finally in Section \ref{sec6} we provide two examples, in one of them the curvature of the manifold in consideration  is identically zero and the other one  is nonpositive but not identically zero.  In both examples,  the proximal point method for  nonconvex problems, proposal in this paper, is applied.

\subsection{Notation and terminology} \label{sec2}
In this section, we introduce some fundamental  properties and notations about Riemannian geometry. These basics facts can be found in any introductory book on Riemannian geometry, such as \cite{MP92} or \cite{S96}.

Let $M$ be a $n$-dimentional connected manifold. We denote by $T_pM$ the $n$-dimentional {\it tangent space} of $M$ at $p$, by $TM=\cup_{p\in M}T_pM$ {\itshape{tangent bundle}} of $M$ and by ${\cal X}(M)$ the space of smooth vector fields over $M$. When $M$ is endowed with a Riemannian metric $\langle \,,\, \rangle$, with the corresponding norm denoted by $\| \; \|$, then $M$ is now a Riemannian manifold. Recall that the metric can be used to define the length of piecewise smooth curves $\gamma:[a,b]\rightarrow M$ joining $p$ to $q$, i.e., such that $\gamma(a)=p$ and $\gamma(b)=q$, by:
\[
l(\gamma)=\int_a^b\|\gamma^{\prime}(t)\|dt,
\]
and, moreover, by minimizing this length functional over the set of all such curves, we obtain a Riemannian distance $d(p,q)$ inducing the original topology on $M$. The metric induces a map $f\mapsto\grad f\in{\cal X}(M)$ which, for each function smooth over $M$, associates its gradient via the rule $\langle\grad f,X\rangle=d f(X),\ X\in{\cal X}(M)$. Let $\nabla$ be the Levi-Civita connection associated with $(M,{\langle} \,,\, {\rangle})$. In each point $p\in M$, we have a linear map $A_X(p) \colon T_pM \to T_pM$ defined by:
\begin{equation}
A_X(p)v=\nabla_vX.
\end{equation}
If $X=\grad f$, where $f\colon M \to\mathbb{R}$ is a twice differentiable function, then $A_X(p)$ is the  {\it Hessian\/} of $f$ at $p$ and is denoted by $\text{Hess }_p f$. A vector field $V$ along $\gamma$ is said to be {\it parallel} if $\nabla_{\gamma^{\prime}} V=0$. If $\gamma^{\prime}$ itself is parallel we say that $\gamma$ is a {\it geodesic}. Given that the geodesic equation $\nabla_{\ \gamma^{\prime}} \gamma^{\prime}=0$ is a second order nonlinear ordinary differential equation, we conclude that the geodesic $\gamma=\gamma _{v}(.,p)$ is determined by its position $p$ and velocity $v$ at $p$. It is easy to check that $\|\gamma ^{\prime}\|$ is constant. We say that $ \gamma $ is {\it normalized} if $\| \gamma ^{\prime}\|=1$. The restriction of a geodesic to a  closed bounded interval is called a {\it geodesic segment}. A geodesic segment joining $p$ to $q$ in $ M$ is said to be {\it minimal} if its length is equals $d(p,q)$ and the geodesic in question is said to be a {\it minimizing geodesic}. If $\gamma$ is a geodesic joining points $p$ and $q$ in $ M$ then, for each $t\in [a,b]$, $\nabla$ induces a linear isometry, relative to ${ \langle}\, ,\, {\rangle}$, $P_{\gamma(a)\gamma(t)}:T_{\gamma(a)}M\to T_{\gamma(t)}M$, the so-called {\it parallel transport} along $\gamma$ from $\gamma(a)$ to $\gamma(t)$. The inverse map of $P_{\gamma(a)\gamma(t)}$ is denoted by $P_{\gamma(a)\gamma(t)}^{-1}:T_{\gamma(t)} M \to T_{\gamma(a)}M$. In the particular case of $\gamma$ to be the unique geodesic segment joining $p$ and $q$, then the parallel transport along $\gamma$ from $p$ to $q$ is denoted by $P_{pq}:T_{p}M\to T_{q}M$.

A Riemannian manifold is {\it complete} if the geodesics are defined for any
values of $t$. Hopf-Rinow's theorem asserts that if this is the case then
any pair of points, say $p$ and $q$, in $M$ can be joined by a (not
necessarily unique) minimal geodesic segment. Moreover, $( M, d)$ is a
complete metric space so that bounded and closed subsets are compact. From the completeness of the Riemannian manifold $M$, the {\it exponential map} $exp_{p}:T_{p}  M \to M $ is defined by $exp_{p}v\,=\, \gamma _{v}(1,p)$, for each $p\in M$.

We denote by $R$ {\it the curvature tensor \/} defined by
$R(X,Y)=\nabla_{X}\nabla_{Y}Z- \nabla_{Y}\nabla_{X}Z-\nabla_{[Y,X]}Z$, with $X,Y,Z\in{\cal X}(M)$, where $[X,Y]=YX-XY$. Moreover, the {\it sectional curvature \/} with respect
to $X$ and $Y$ is given by $K(X,Y)=\langle R(X,Y)Y , X\rangle /(||X||^{2}||X||^{2}-\langle X\,,\,Y\rangle ^{2})$, where $||X||=\langle X,X\rangle ^{2}$. If $K(X,Y)\leqslant 0$ for all $X$ and $Y$, then $M$ is called  a {\it Riemannian manifold of nonpositive curvature \/} and we use the short notation $K\leqslant 0$.
\begin{theorem} \label{T:Hadamard}
Let $M$ be a complete, simply connected Riemannian manifold with nonpositive sectional
curvature. Then $M$ is diffeomorphic to the Euclidean space $\mathbb{R}^n $, $ n=dim M $. More
precisely, at any point $p\in M $, the exponential map $exp_{p}$ is a diffeomorphism.
\end{theorem}
\begin{proof}
See Lemma 3.2 of \cite{MP92}, p. 149 or Theorem 4.1 of \cite{S96}, p. 221.
\end{proof}

 A complete simply connected Riemannian manifold of nonpositive sectional curvature is called a {\it{Hadamard manifold}}. Thus Theorem \ref{T:Hadamard} states that if $M$ is a Hadamard manifold, then $M$ has the same topology and differential structure of the Euclidean space $\mathbb{R}^n$. Furthermore, are known some similar geometrical properties of the Euclidean space $\mathbb{R}^n$, such as, given two points there
exists an unique geodesic segment that joins them. {\it In this paper, all manifolds $M$ are assumed to be Hadamard and finite dimensional}.

\section{Convexity in Hadamard manifolds} \label{sec3}
In this section, we introduce some fundamental properties and notations of convex analysis on Hadamard manifolds. References of the convex analysis, on the Euclidean space $\mathbb{R}^{n}$ may be found in \cite{HL93} and on Riemannian manifolds may be found in \cite{XFL2002}, \cite{FO2000}, \cite{RAP97}, \cite{S96}, \cite{S94} and  \cite{U94}.

The set $\Omega\subset M$ is said to be {\it convex \/} if for any geodesic segment, with end points in $\Omega$, is contained in $\Omega$. Let $\Omega\subset M$  be an open convex set. A function $f:M\to\mathbb{R}$ is said to be {\it convex} (respectively, {\it strictly convex}) on $\Omega$ if for any geodesic segment $\gamma:[a, b]\to\Omega$ the composition $f\circ\gamma:[a, b]\to\mathbb{R}$ is convex (respectively, strictly convex). Moreover, a function $f:M\to\mathbb{R}$ is said to be {\it strongly convex} on $\Omega$ with constant $L>0$ if, for any geodesic segment $\gamma:[a, b]\to\Omega$, the composition $f\circ\gamma:[a, b]\to\mathbb{R}$ is strongly convex with constant $L\|\gamma'(0)\|^2$. Take $p\in M$. A vector $s \in T_pM$ is said to be a {\it subgradient\/} of $f$ at $p$, if:
\[
f(q) \geq f(p) + \langle s,\exp^{-1}_pq\rangle,
\]
for any $q\in M$. The set of all subgradients of $f$ at $p$, $\partial f(p)$, is called the {\it subdifferential\/} of $f$ at $p$.

The following result provides a characterization of convexity in the case of differentiable functions.
\begin{proposition}\label{PropoFC1}
Let $\Omega\subset M$ be an open convex set and $f:M\to\mathbb{R}$ a differentiable function on $\Omega$. We say that $f$ is convex on $\Omega$ if, and only if, for any $p\in\Omega$:
\[
f(q)-f(p)\geq \langle \grad f(p), \exp^{-1}_pq\rangle,
\]
for all $q\in\Omega$.
\end{proposition}
\begin{proof}\ See Theorem 5.1 of \cite{U94}, page 78.
\end{proof}

The most important consequence of the previous proposition is that with $f$ being convex, any of its critical points are global minimum points. In particular, if $M$ is compact, then $f$ is constant. Moreover, $0\in\partial f(p)$ if, and only if, $p$ is a minimum point of $f$ in $M$. See, for example, \cite{U94}.
\begin{definition} \label{def4.2}
Let $\Omega\subset M$ be an open convex set and $X$ a vector field defined in $M$. $X$ is said to be {\it monotone\/} on $\Omega$, if:
\begin{equation}\label{eq4.2}
\left\langle \exp_{q}^{-1}p\, , \,P_{qp}^{-1} X(p)-X(q)\right\rangle \geq 0, \qquad p,\,q\in\Omega,
\end{equation}
where $P_{qp}$ is the parallel transport along the geodesic  joining  $q$ to $p$. If (\ref{eq4.2}) is satisfied with strict inequality for all $p,q\in\Omega$, $p\neq q$, then $X$ is said to be {\it strictly monotone\/}. Moreover,   $X$ is strongly monotone if there exists  $\lambda>0$ such that:
\begin{equation}\label{eq4.3}
\left\langle \exp_{q}^{-1}p\,,\,P_{qp}^{-1} X(p)-X(q)\right\rangle \geq \lambda d^2(p,q)\qquad p,\,q\in\Omega.
\end{equation}
\end{definition}

\begin{remark}
In the particular case that $M=\mathbb{R}^n$ with the usual metric, inequality \eqref{eq4.2} and \eqref{eq4.3} becomes, respectively:
\[
\left\langle p-q\,,\,X(p)-X(q)\right\rangle \geq 0,  \qquad  \left\langle p-q\,,\,X(p)-X(q)\right\rangle \geq\lambda\|p-q\|^2,
\]
because $\exp_{q}^{-1}p=p-q$ and $P_{qp}^{-1}=I$. Therefore the Definition \ref{def4.2} extends the concept of monotone operators from $\mathbb{R}^n$ to Riemannian manifolds.
\end{remark}

Now we present an important example of strong monotone vector field being particularly useful in the remainder of this work.

Take ${p}\in M $ and let $exp^{-1}_{p}:M\to T_{p}M$ be the inverse of the exponential map. Note that $d({q}\, , \, p)\,=\,||exp^{-1}_{p}q||$,  the map $d^2(\,.\,,{p})\colon M\to\mathbb{R}$
is  $C^{\infty}$ and
\[
\grad \frac{1}{2}d^2(q,{p})=-exp^{-1}_{q}{p},
\]
($M$ is a Hadamard manifold). See, for example, Proposition 4.8 of \cite{S96}, p. 108.
\begin{proposition} \label{P:dg}
Take  ${p}\in M$. The gradient vector field $\grad (d^2(\,.\,,{p})/2)$ is strongly monotone with $\lambda=~1$.
\end{proposition}
\begin{proof}
See \cite{XFL2002}.
\end{proof}
\begin{proposition} \label{prop:fconv}
Let $\Omega\subset M$ be an open convex set and $f:M\to\mathbb{R}$ a differentiable function on $\Omega$.
\begin{itemize}
\item[(i)] $f$ is convex on $\Omega$ if and only if $\grad f$ is monotone on $\Omega$;
\item[(ii)] $f$ is strictly convex on $\Omega$ if  and only if $\grad f$ is strictly monotone on $\Omega$;
\item[(iii)] $f$ is strongly convex on $\Omega$ if  and only if $\grad f$ is strongly monotone on $\Omega$.
\end{itemize}
\end{proposition}
\begin{proof} See \cite{XFL2002}.
\end{proof}
\begin{remark}
Take ${p}\in M$. From Propositions \ref{P:dg} and \ref{prop:fconv} it follows that the map $d^2(\,.\,,{p})/2$ is strongly convex.
\end{remark}

\begin{proposition} \label{prop:supfcov}
Let $\Omega\subset M$ be a convex set and $T\subset \mathbb{R}$ a compact set.
Let $\psi:M\times T\to \mathbb{R}$ be a continuous function on $\Omega\times T$ such that $\psi_{\tau}:=\psi(.,\tau):M\to \mathbb{R}$ is strongly convex on $\Omega$ with constant $L>0$ for all $\tau\in T$. Then, $\phi:M\to\mathbb{R}$ defined by:
\[
\phi(p):=\max_{\tau\in T}\psi(p,\tau),
\]
is strongly convex on $\Omega$ with constant $L$. In particular, if $\psi_{\tau}$ is convex for all $\tau\in T$ then $\phi$ is convex on $\Omega$.
\end{proposition}
\begin{proof}
Since $T$ is compact and $\psi$ is continuous, the function $\phi$ is well defined. Let $\gamma:[a, b]\to\Omega$  be a geodesic segment. Because $\psi_{\tau}$ is strongly convex with constant $L$ for each $\tau\in T$, we have:
\[
(\psi_\tau\circ\gamma)(\alpha t_1+(1-\alpha)t_2)\leq\alpha(\psi_\tau\circ\gamma)(t_1)+(1-\alpha)(\psi_\tau\circ\gamma)(t_2)-\frac{1}{2}(1-\alpha)\alpha L\|\gamma'(0)\|^2,
\]
for all $t_1, t_2 \in [a, b]$ and $\alpha\in [0,1]$. Thus, taking the maximum in $\tau$ in both sides of the above inequality, we obtain:
\[
(\phi\circ\gamma)(\alpha t_1+(1-\alpha)t_2)\leq\alpha(\phi\circ\gamma)(t_1)+(1-\alpha)(\phi\circ\gamma)(t_2)-\frac{1}{2}(1-\alpha)\alpha L\|\gamma'(0)\|^2,
\]
which implies that $\phi\circ\gamma: [a, b]  \to  \mathbb{R}$ is strongly convex with constant $L\|\gamma'(0)\|^2$. So, $\phi$ is strongly convex on $\Omega$ with constant $L$. The proof of the second part is immediate.
\end{proof}

\begin{definition}
Let $\Omega\subset M$ be an open convex set. A function $f: M \to \mathbb{R}$ is said to be Lipschitz on $\Omega$ if there exists a constant $L:=L(\Omega)\geq 0$ such that
\begin{equation}\label{Lipsch1}
|f(p)-f(q)|\leq Ld(p,q), \qquad p,q\in\Omega.
\end{equation}
Moreover, if for each $p_0\in \Omega$ there exists $L(p_0)\geq 0$ and $\delta=\delta(p_0)>0$ such that inequality (\ref{Lipsch1}) holds with $L=L(p_0)$ for all $p,q\in B_{\delta}(p_0):=\{p\in \Omega : d(p,p_0)<\delta\}$, then $f$ is called locally Lipschitz on $\Omega$.
\end{definition}
\begin{definition}
Let $\Omega\subset M$ be an open convex set and $f:M\to\mathbb{R}$ a continuously differentiable function on $\Omega$.  The gradient vector field of $f$, $\grad f$, is said to be Lipschitz with constant $\Gamma \geq 0$ on $\Omega$ whenever:
\[
\|\grad f(q)-P_{pq}\grad f(p)\|\leq\Gamma
d(p,q), \qquad p, q \in \Omega,
\]
where  $P_{pq}$ is the parallel transport along the geodesic joining $p$ to $q$.
\end{definition}
\begin{proposition}\label{gradLips:01}
Let $\Omega\subset M$ be an open convex set and $f:M\to\mathbb{R}$ a twice continuously differentiable function on $\Omega$. If $\text{Hess }_p f$ is bounded on $\Omega$, then the gradient vector field $\grad f$ is Lipschitz on $\Omega$.
\end{proposition}
\begin{proof}
The proof is an immediate consequence of the fundamental theorem of calculus for vector fields, see for example \cite{FS02}.
\end{proof}
\section{Directional derivatives}\label{sec:dd}
In this section we present some properties of the directional derivative of a convex function defined on a Hadamard manifold, including a characterization of the directional derivative and of the subdifferential of the maximum of a certain class of convex functions. We also give a definition of the generalized directional derivative of a locally Lipschitz function (not necessarily convex), see Azagra et.all \cite{Azagra2005}, and an important property of the  subdifferential of the maximum of continuously differentiable functions.
\subsection{Directional derivatives of convex functions}\label{sec:dd1}
In this subsection we present the definition of the directional derivative of a convex function defined on a Hadamard manifold and some properties involving its subdifferential, which  allow us to obtain an important property of the  subdifferential of the maximum of a certain class of convex functions.

Let $\Omega \subset M$ be an open convex set and $f\colon M\to R$  a convex function on $\Omega$. Take $p \in \Omega$,  $v \in T_pM $ and $\delta >0$ and let $\gamma :[-\delta\,, \,\delta]\to \Omega$ be the geodesic segment such that  $\gamma(0)=p$ and $\gamma'(0)=v$. Due to the convexity of $f\circ\gamma:[-\delta\,, \,\delta]\to\mathbb{R}$, the function $q_{\gamma}:(0\,, \,\delta]\to\mathbb{R}$, given by
\begin{equation}\label{eq3.1}
q_{\gamma}(t):=\frac{f(\gamma(t)) -f(p) }t,
\end{equation}
is nondecreasing. Moreover, since $f$ is locally Lipschitzian, it follows that $q_{\gamma}$ is bounded near
zero. This leads to the following definition:
\begin{definition}
Let $\Omega \subset M$ be an open convex set and $f\colon M\to R$  a convex function on $\Omega$. Then the {\it directional derivative} of $f$ at $p\in \Omega$ in the direction of $v\in T_{p}M$ is defined by
\begin{equation} \label{eq3.2}
f'(p,v):= \lim_{t \to 0^+} q_{\gamma}(t)=\inf_{t>0}q_{\gamma}(t),
\end{equation}
where $\delta>0$ and   $\gamma :[-\delta\,, \,\delta]\to \Omega$  is the geodesic segment  such that $\gamma(0)=p$ and $\gamma'(0)=v$.
\end{definition}
\begin{proposition}\label{thm3.9}
Let $\Omega\subset M$ be an open convex set and $f:M\to\mathbb{R}$ a convex function on $\Omega$. Then, for each fixed $p\in\Omega$, the subdifferential $\partial f(p)$ is convex.
\end{proposition}
\begin{proof}
See Theorem 4.6 of \cite{U94}, p. 74.
\end{proof}
\begin{proposition}\label{prop3.11}
Let $\Omega\subset M$ be an open convex set and $f:M\to\mathbb{R}$ a convex function on $\Omega$. Then, for each point fixed $p\in\Omega$, the following statement holds:
\begin{itemize}
\item[i)] $f'(p,v)=\max_{s \in \partial f(p)} \langle s,v\rangle$, for all $v \in T_{p}M$;
\item[ii)] $\partial f(p)=\left\{s \in T_{p}M : f'(p,v) \geq \langle s,v \rangle,  v \in
T_{p}M \right\}$.
\end{itemize}
\end{proposition}
\begin{proof}
See \cite{XFL2000}.
\end{proof}
\begin{proposition} \label{p:env}
Let $T$ be a compact set, $\Omega\subset M$ an open convex set and $h:M\times T\to \mathbb{R}$ a continuous function on $\Omega\times T$ such that $h(.\,,\tau):M\to \mathbb{R}$ is convex on $\Omega$ for all $\tau\in T$. If $f:M\to \mathbb{R}$ is given by $f(p)=\max_{\tau\in T} h(p,\tau)$, then $f$ is convex on $\Omega$ and
\[
f'(p,v)=\max  \limits_{\tau\in T(p)}h'(p,\tau,v), \qquad \  p\in\Omega, \quad  v\in T_pM,
\]
where $T(p)=\{\tau\in T : f(p)=h(p,\tau)\}$. Moreover, if $h(.,\tau)$ is differentiable on $\Omega$ for all $\tau\in T$ and $\grad_ph(p,.)$ is continuous for all $p\in\Omega$, then:
\[
\partial f(p)=\conv \left\{\grad_p h(p,\tau): \tau\in T(p)\right\}.
\]
\end{proposition}
\begin{proof}
Since $T$ is compact  $f$ is well defined and  its convexity  follows from Proposition~\ref{prop:supfcov}. Now, take $p\in\Omega$,  $v\in T_pM$ and the geodesic segment $\gamma:[-\delta,\delta]\to\Omega$, $\delta>0$,  such that $\gamma(0)=p$ and $\gamma'(0)=v$. Using that $T$ is compact we have $T(p)\neq\emptyset$. Hence, taking $\tau\in T(p)$ we obtain from definition of $f$ and $T(p)$
\[
\frac{f(\gamma(t))-f(p)}t\geq\frac{h(\gamma(t),\tau)-h(p,\tau)}t, \qquad t\in (0,\delta).
\]
Given that $f$ and $h(.,\tau)$ are convex, letting $t$ goes to $0$, the above inequality yields:
\[
f'(p,v)\geq h'(p,\tau,v),\qquad  p\in\Omega,  \;   v\in T_pM,\ \;\tau\in T(p).
\]
Therefore,
\begin{equation}\label{Desi110}
f'(p,v)\geq\sup\limits_{\tau\in T(p)}h'(p,\tau,v), \qquad  p\in\Omega,  \; v\in T_pM.
\end{equation}
Now, we are going to prove the equality in the above equation. Let $\{t_k\}\subset (0,\delta)$ such that $t_k$ converges to $0$ as  $k$ goes to $+\infty$. We define
\begin{equation} \label{eq:gtk}
p^k:=\gamma(t_k), \qquad  \tau_k\in T(p^k).
\end{equation}
The last equality implies that $f(\gamma(t_k))=h(\gamma(t_k),\tau_k)$. Therefore, as $f$ is convex, combining \eqref{eq3.2} and the definition of $f$, we obtain:
\[
f'(p,v)\leq\frac{f(\gamma(t_k))-f(p)}{t_k}\leq \frac{h(\gamma(t_k),\tau_k)-h(p,\tau_k)}{t_k}.
\]
Since $\{\tau_k\}\subset T$ and $T$ is compact, we can suppose (taking a subsequence, if necessary) that it converges to $\bar{\tau}\in T$ as  $k$ goes to $+\infty$. Thus, letting  $k$ goes to $+\infty$ in the latter inequality, we have:
\[
f'(p,v)\leq h'(p,\bar{\tau},v),
\]
Since $h$ is continuous and $x \mapsto h (x , \tau)$  is convex for each $\tau\in T$. Note that, if $\bar{\tau}\in T(p)$, then the last inequality implies that \eqref{Desi110} holds with equality. So, in order to prove the first part it is sufficient to prove that $\bar{\tau}\in T(p)$. First note that using \eqref{eq:gtk} and the definitions of  $T(p^k)$ and $f$, we conclude that:
\[
h(p^k,\tau_k)=\max\limits_{\tau\in T}h(p^k,\tau)\geq h(p^k,\tau),\qquad  \tau\in T.
\]
Hence, letting $k$ goes to $+\infty$, we have $h(p,\bar{\tau})\geq h(p,\tau)$, for all $\tau\in T$, which, combined with the definition of $f$, give $f(p)=h(p,\bar{\tau})$, and the first part is completed.

To  prove the second part, take $p\in\Omega$ and $\tau\in T(p)$. From Proposition \ref{PropoFC1}, the convexity of $h(.,\tau)$ implies:
\[
h(q,\tau)\geq h(p,\tau)+\langle\grad_p h(p,\tau),\exp^{-1}_pq\rangle,\qquad q\in\Omega.
\]
Because $\tau\in T(p)$ we have $h(p,\tau)=f(p)$, which along with the definition of $f$ as well as with latter equation yields:
\[
f(q)\geq f(p)+\langle\grad_p h(p,\tau),\exp^{-1}_pq\rangle.
\]
So,  $\grad_p h(p,\tau)\in\partial f(p)$. Since Proposition \ref{thm3.9} implies that $\partial f(p)$ is convex,  we conclude that:
\[
\conv\{\grad_p h(p,\tau): \tau\in T(p)\}\subseteq\partial f(p).
\]
We claim that this inclusion holds with equality. Indeed, assume by contradiction that:
\[
\exists\ y\in\partial f(p), \quad  y\notin\conv\{\grad_p h(p,\tau): \tau\in T(p)\}.
\]
Due the fact that  $\grad_p h(p,.)$ is continuous and $T(p)$ is a compact set we conclude that the set  $\conv\{\grad_p h(p,\tau): \tau\in T(p)\}$  is  compact. Thus, by the Separation Theorem on $T_pM$, there exists $v\in T_pM-\{0\}$ and $a\in\mathbb{R}$ such that
\[
\langle y,v\rangle>a>\langle\grad_p h(p,\tau), v\rangle, \quad\forall\, \tau\in T(p).
\]
Since $h'(p,\tau,v)=\langle\grad_p h(p,\tau), v)\rangle$, it follows from the latter inequality and first part of the proposition that
\[
\langle y, v\rangle>\max\limits_{\tau\in T(p)}h'(p,\tau, v)=f'(p, v).
\]
Since $y\in\partial f(p)$, we obtain a contradiction with Proposition \ref{prop3.11} i. This proves our claim and conclude the proof.
\end{proof}
\begin{corollary}\label{CoroMPP1}
Let $\Omega\subset M$ be a open convex set and $h_i:M\to\mathbb{R}$ a differentiable convex function on $\Omega$ for $i\in I:=\{ 1,...,m\}$. If $h:M\to \mathbb{R}$ is  defined by $h(p):=\max_{i\in I}h_i(p)$, then:
\[
\partial h(p)=\conv\{\grad h_i : i\in I(p)\}=\left\{y\in T_pM :  y=\sum\limits_{i\in I(p)}\alpha_i\grad h_i(p), \sum\limits_{i\in I(p)}\alpha_i=1, \alpha_i\geq 0\right\},
\]
where $I(p):=\{i\in I: h(p)=h_i(p)\}$. In particular, $p$ minimizes $h$ on $\Omega$, if and only if, there exist $\alpha_i\geq 0$, $i\in I(p)$, such that:
\[
0=\sum\limits_{i\in I(p)}\alpha_i\grad h_i(p), \qquad \sum\limits_{i\in I(p)}\alpha_i=1.
\]
\end{corollary}
\begin{proof}
It follows directly from Proposition \ref{p:env}.
\end{proof}
\subsection{Directional derivatives of locally Lipschitz functions}\label{sec:dd2}
In the sequel we present the definition of  generalized directional derivative of a locally Lipschitz function (not necessarily convex) and an important property of the  subdifferential of the maximum of continuously differentiable functions.
\begin{definition}
Let $\Omega\subset M$ be an open convex set and $f:M\to\mathbb{R}$ a locally Lipschitz function on $\Omega$. The generalized directional derivative of $f$ at $p\in \Omega$ in the direction $v\in T_pM$ is defined by:
\begin{equation}\label{Clarke1}
f^\circ(p,v):=\limsup\limits_{t\downarrow 0\ q\to p}\frac{f(\exp_qt(D\exp_p)_{\exp^{-1}_pq}v))-f(q)}{t},
\end{equation}
where $(D\exp_p)_{\exp^{-1}_pq}$ denotes the differential of $\exp_p$ at $\exp^{-1}_pq$.
\end{definition}
It is worth noting that an equivalent definition has appeared in \cite{Azagra2005}.
\begin{remark}
The generalized directional derivative is well defined. Indeed,
let $L_p>0$ be the Lipschitz  constant of $f$ in $p$ and $\delta=\delta(p)>0$
such that
\[
|f(\exp_q t(D\exp_p)_{\exp^{-1}_pq}v)-f(q)|\leq L_p \,d(\exp_q t(D\exp_p)_{\exp^{-1}_pq}v, \,q), \quad q\in B_{\delta}(p), \quad t\in[0,\delta).
\]
Since $d(\exp_q t(D\exp_p)_{\exp^{-1}_pq}v, \,q)= t \| (D\exp_p)_{\exp^{-1}_pq}v\|$, the above inequality becomes:
\[
|f(\exp_q t(D\exp_p)_{\exp^{-1}_pq}v)-f(q)|\leq L_p \, t \| (D\exp_p)_{\exp^{-1}_pq}v\|, \quad q\in B_{\delta}(p), \quad t\in[0,\delta).
\]
Since $\lim_{ q\to p}\,(D\exp_p)_{\exp^{-1}_pq}v=v$, our statement follows from the latter inequality.
\end{remark}

\begin{remark}
Note that, if $M=\mathbb{R}^n$ then  $\exp_{p}w=p+w$ and
\[
D(\exp_p)_{\exp^{-1}_pq}v=v.
\]
In this case, \eqref{Clarke1} becomes:
\[
f_E^{\circ}(p,v)=\limsup\limits_{t\downarrow 0\ q\to p}\frac{f(q+tv)-f(q)}t,
\]
which is the Clarke's generalized directional derivative in the Euclidean case, see \cite{Clarke1983}. Therefore, the generalized differential derivative on Hadamard manifolds is a natural extension of the Clarke's generalized differential derivative.
\end{remark}
Next we generalize the definition of subdifferential for locally Lipschitz functions defined on Hadamard manifolds, see Proposition \ref{prop3.11} item ii.
\begin{definition}
Let $\Omega\subset M$ be an open convex set and $f:M\to\mathbb{R}$ a locally Lipschitz function on $\Omega$. The generalized subdifferential of $f$ at $p\in \Omega$, denoted by $\partial^{\circ} f(p)$, is defined by:
\[
\partial^{\circ}f(p):=\{w\in T_pM : f^{\circ}(p,v)\geq\langle w,v\rangle \mbox{ for all}\  v\in T_pM\}.
\]
\end{definition}

\begin{remark}\label{Remark:Azagra2005}
If $f$ is convex on $\Omega$, then $f^{\circ}(p,v)=f'(p,v)$ (respectively, $\partial^\circ f(p)=\partial f(p)$) for all $p\in \Omega$, i.e., the directional derivatives (respectively,  subdifferential) for Lipschitz functions is a generalization of the directional derivatives (respectively,  subdifferential) for convex functions.
See  \cite{Azagra2005} Claim $5.4$ in the proof of Theorem $5.3$.
\end{remark}

\begin{definition}
Let $\Omega\subset M$ be an open convex set and $f:M\to\mathbb{R}$ a locally Lipschitz function on $\Omega$. A point $p\in \Omega$ is a stationary point of $f$ if $0\in \partial^{\circ}f(p)$.
\end{definition}
\begin{lemma}\label{lemacl1}
Let $\Omega\subset M$ be an open set. If $f:M\to\mathbb{R}$ is locally Lipschitz on $\Omega$ and $g:M\to\mathbb{R}$ is continuously differentiable on $\Omega$, then:
\begin{equation}\label{ddg:1}
(f+g)^{\circ}(p,v)=f^{\circ}(p,v)+g'(p,v)\qquad p\in \Omega, \quad v\in T_pM.
\end{equation}
As a consequence,
\begin{equation}\label{sdg:1}
\partial^{\circ}(f+g)(p)=\partial^{\circ}f(p)+ \grad g(p),\qquad p\in\Omega.
\end{equation}
\end{lemma}
\begin{proof}
Using the definition of the generalized directional derivative and simple algebraic manipulations, we obtain:
\[
(f+g)^{\circ}(p,v)=\limsup\limits_{t\downarrow 0\ q\to p}
\left[\frac{f(\exp_qt(D\exp_p)_{\exp^{-1}_pq}v)-f(q)}t+\frac{g(\exp_qt(D\exp_p)_{\exp^{-1}_pq}v)-g(q)}t\right].
\]
From  upper limit properties together with the definition of the generalized directional derivative and Remark \ref{Remark:Azagra2005}, we obtain
\begin{equation}\label{P:ddg1}
(f+g)^{\circ}(p,v)\leq f^{\circ}(p,v)+g'(p,v).
\end{equation}
Because $f^{\circ}(p,v)=\left((f+g)+(-g)\right)^{\circ}(p,v)$, above inequality implies in particular  that
\[
f^{\circ}(p,v)\leq (f+g)^{\circ}(p,v)+(-g)'(p,v),
\]
which is equivalent to
\[
(f+g)^{\circ}(p,v)\geq f^{\circ}(p,v)+g'(p,v).
\]
So, last inequality together  with inequality \eqref{P:ddg1} yields the equality \eqref{ddg:1}.

The proof of the equality \eqref{sdg:1} is an immediate consequence of the equality \eqref{ddg:1} and the definition of the generalized subdifferential.
\end{proof}
It is possible to prove that the next result holds with equality. However we will prove just the inclusion needed to prove our main result.
\begin{proposition}\label{SubClarke1}
Let $\Omega\subset M$ be an open convex set and   $I=\{1,...,m\}$.
Let $f_{i}:M\to \mathbb{R}$ be a continuously differentiable function on $\Omega$ for all $i\in I$ and $f:M \to\mathbb{R}$  defined by
\[
f(p):=\max_{i\in I} f_i(p).
\]
Then $f$ is Lipschitz locally on $\Omega$ and for each $p\in \Omega$
\[
\conv\{\grad f_i(p) : i\in I(p)\}\subset\partial^{\circ} f(p),
\]
where $I(p):=\{i\in I : f_i(p)=f(p)\}$.
\end{proposition}
\begin{proof}
Since $f_i$ is continuously differentiable on $\Omega$ we conclude that
$f_i$ is Lipschitz locally in $\Omega$, for all $i\in I$. Thus, for each $\tilde{p}\in \Omega$ and $i\in I$,  there exists $\delta_i, L_i>0$ such that
\[
|f_i(p)-f_i(q)|\leq L_id(p,q), \qquad p\,,q\in B_{\delta_i}(\tilde{p}).
\]
On the other hand,
\[
|\max_{i\in I}f_i(p)-\max_{i\in I}f_i(q)|\leq\max_{i\in I}|f_i(p)-f_i(q)|.
\]
Combining two last equations with the definition of $f$, we obtain:
\[
|f(p)-f(q)|\leq Ld(p, q) \qquad p,q\in B_{\delta}(\tilde{p}),
\]
where $\delta=\min_{i\in I}\delta_i$ and $L=\max_{i\in I}L_i$, which prove the first part.

In order to prove the second part, take $p\in \Omega$, $u\in\conv\{\grad f_i(p) : i\in I(p)\}$ and $v\in T_{p}M$. Then, there exist a constant $\alpha_i\geq 0$ for $i\in I(p)$ with $\sum_{i\in I(p)}\alpha_i=1$ such that
\[
u=\sum_{i\in I(p)}\alpha_i\grad f_i(p).
\]
Since $f_i$ is differentiable for all $i\in I$, simple algebraic manipulation yields
\[
\langle u,v\rangle=\sum_{i\in I(p)}\alpha_i\langle\grad f_i(p),v\rangle=\sum_{i\in I(p)}\alpha_if'_i(p, v).
\]
As $f$ is locally Lipschitz at $p$, the definitions of $f$, $I(p)$ and generalized directional derivative imply:
\[
f'_i(p, v)\leq f^{\circ}(p, v),
\]
which, together with the latter equation, gives $\langle u,v\rangle\leq f^{\circ}(p, v)$, and the proof follows from the definition of $\partial^{\circ}f(p)$.
\end{proof}

\section{Proximal Point Method for Nonconvex Problems}\label{sec5}
In this section we present an application of the proximal point method to minimize a real-valued function (not necessarily convex) given by the maximum of a certain class of  continuously differentiable functions. Our goal is to prove the following theorem:
\begin{theorem}\label{MPP10}
Let $\Omega\subset M$ be an open convex set,  $q\in M$ and   $I=\{1,...,m\}$.
Let $f_{i}:M\to \mathbb{R}$ be a continuously differentiable function on $\Omega$ and continuous on $\bar{\Omega}$ (closure of $\Omega$), for all $i\in I$, and $f:M\to\mathbb{R}$ defined by
\[
f(p):=\max_{i\in I} f_i(p).
\]
Assume that  $-\infty<\inf_{p\in M}f(p)$, $\grad f_i$ is Lipschitz on $\Omega$ with constant $L_{i}$ for each $i\in I$ and
\[
L_f(f(q))=\left\{p\in M: f(p)\leq f(q)\right\}\subset \Omega,\qquad \inf_{p\in M}f(p)<f(q).
\]
Take $0<\bar{\lambda}$ and a sequence $\{ \lambda_{k}\}$ satisfying $\max\limits_{i\in I}L_i<\lambda_k\leq \bar{\lambda}$ and $\hat{p}\in L_f(f(q))$. Then the proximal point method
\begin{equation}\label{E:1.22}
p^{k+1}:=\argmin_{p\in M} \left\{f(p)+\frac{\lambda_k}{2}d^2(p,{p^k})\right\}, \qquad k=0, 1, \ldots,
\end{equation}
with starting point $p^0=\hat{p}$ is well defined, the generated sequence $\{p^{k}\}$ rest in $L_f(f(q))$ and satisfies only one of the following statements
\begin{itemize}
\item[i)] $\{p^{k}\}$ is finite, i.e., $p^{k+1}=p^k$ for some $k$ and,  in this case, $p^k$ is a stationary point of $f$,
\item[ii)] $\{p^{k}\}$ is infinite and, in this case, any accumulation point of $\{p^k\}$ is a stationary point of $f$.
\end{itemize}
Moreover, assume that  the minimizer set of $f$ is non-empty, i. e.,
\begin{itemize}
\item[{\bf h1)}] $U^*=\{p : f(p)=\inf_{p\in M}f(p)\}\neq\emptyset$.
\end{itemize}
Let $c\in(\inf_{p\in M}f(p),\; f(q))$. If, in addition,  the following assumptions hold:
\begin{itemize}
\item[{\bf h2)}] $L_f(c)$ is convex and  $f$ is convex on $L_f(c)$;
\item[{\bf h3)}] For all $p\in L_f(f(q))\setminus L_f(c)$ and $y(p)\in\partial^{\circ}f(p)$ we have $\|y(p)\|>\delta>0, $
\end{itemize}
then the sequence $\{p^{k}\}$ generated by \eqref{E:1.22} with
\begin{equation}\label{mpp102}
\max\limits_{i\in I}L_i<\lambda_k\leq\bar{\lambda}, \qquad k=0, 1, \ldots
\end{equation}
converge to a  point $p^*\in U^*$.
\end{theorem}
\begin{remark}
The continuity of each function $f_i$ on $\bar{\Omega}$ in {\bf h2} guarantees that the level sets of the fuction $f$, in particular the solution set $U^*$, are closed in the topology of the manifold.
\end{remark}

In the next remark we show that if $\Omega$ is bounded and $f_i$ is convex on $\Omega$,  for all $i\in I$, then $f$ satisfies the assumptions {\bf h2} and {\bf h3}.
\begin{remark}\label{re:tconver}
If $f_i$ is also a convex function on $\Omega$ for each $i\in I$ then by the Proposition \ref{prop:supfcov}, the function $f$ is convex on $\Omega$ and the assumption {\bf h2} is satisfied for all $c\leq f(q)$. Moreover, from Remark \ref{Remark:Azagra2005},
\begin{equation} \label{eq:rsfr}
\partial f^{\circ}(p)=\partial f(p),\quad\forall\; p\in \Omega.
\end{equation}
Take $c\in (\inf_{p\in M}f(p),\; f(q))$ and let us suppose that {\bf h1} hold and $\Omega$ is a bounded set. Then, we have
\begin{equation}\label{CconvexMPP1}
0<\sup\left\{d(p^*,p): p^*\in U^*, p\in  L_f(f(q))\setminus L_f(c)\right\}=\epsilon<+\infty.
\end{equation}
Let $p^*\in U^*$ be fixed, $p\in L_f(f(q))\setminus L_f(c)$ and $y(p)\in \partial f(p)$. The convexity of $f$ on $\Omega$ implies that:
\[
\langle y(p)\;,\;-\exp^{-1}_pp^*\rangle\geq f(p)-f(p^*).
\]
Since $\|y(p)\|\|\exp^{-1}_p{p^*}\|\geq\langle y(p), -\exp^{-1}_p{p^*}\rangle$, $d(p^*,p)=\|\exp^{-1}_p{p^*}\|$, $p\in L_f(f(q))\setminus L_f(c)$ and $U^*$ is a proper subset of $L_f(c)$, from the above inequality, we obtain:
\[
\|y(p)\|d(p^*, p)> c-f(p^*)>0.
\]
Thus, from \eqref{CconvexMPP1} and latter inequality
\[
\|y(p)\|\epsilon > c-f(p^*)>0.
\]
Therefore, choosing $\delta = (c-f(p^*))/\epsilon$, we have:
\[
\|y(p)\|>\delta>0,
\]
which, combined with \eqref{eq:rsfr}, shows that $f$ satisfies {\bf h3}.
\end{remark}

In order to prove proving the above theorem we need of some preliminary results. From now on we assume that all assumptions on Theorem \ref{MPP10} holds, with the exception of {\bf h1}, {\bf h2} and {\bf h3}, which will be assumed to hold
only when explicitly stated.
\begin{lemma}\label{MPP8}
For all $\tilde{p}\in M$ and $\lambda$ satisfying
\[
 \sup_{i \in I}L_{i}<\lambda ,
\]
the function $f_i+(\lambda/2) d^2(. \,, \tilde{p} )$ is strongly convex in $\Omega$ with constant $\lambda - \sup_{i \in I}L_{i}$. Consequently, $f+(\lambda/2) d^2(. \,, \tilde{p} )$ is strongly convex on $\Omega$ with constant $\lambda - \sup_{i \in I}L_{i}$.
\end{lemma}
\begin{proof}
Due to the finiteness of $I$, the function $f$ is well defined. Take $i \in I$,  $\tilde{p}\in M$ and define  $h_{i}:=f_i+(\lambda/2) d^2(. \,, \tilde{p} )$. Note that $\grad h_{i}(p)= \grad f_i(p)-\lambda \exp^{-1}_{p}{ \tilde{p}}$. Thus, for all $p, q\in\Omega$:
\begin{multline*}
 \langle P^{-1}_{qp}\grad h_{i}(p)-\grad h_{i}(q), \;\exp^{-1}_{q}p\rangle= \langle P^{-1}_{qp}\grad f_i (p)-\grad f_i(q),\; \exp^{-1}_{q}p\rangle \\ - \lambda \langle P^{-1}_{qp}\exp^{-1}_{p}{\tilde{p}} - \exp^{-1}_{q}{ \tilde{p}}, \; \exp^{-1}_{q}p \rangle.
\end{multline*}
Since $\langle P^{-1}_{qp}\grad f_i (p)-\grad f_i(q), \exp^{-1}_{q}p\rangle \geq -\|P^{-1}_{qp}\grad f_i (p)-\grad f_i(q)\| \|\exp^{-1}_{q}p\|$, using equality $d(p,q)=\|\exp^{-1}_{q}p\|$, Proposition \ref{P:dg} and the above equation, we obtain:
\[
\langle P^{-1}_{qp}\grad h_{i}(p)-\grad h_{i}(q), \exp^{-1}_{q}p\rangle \geq -\|P^{-1}_{qp}\grad f_i (p)-\grad f_i(q)\|d(q,p)+\lambda d^{2}(p,q).
\]
Now, as $\grad f_i$ is Lipschitz on $\Omega$ with constant $L_{i}$ and the parallel transport is an isometry, the latter equation becomes
\[
\langle P^{-1}_{qp}\grad h_{i}(p)-\grad h_{i}(q), \exp^{-1}_{q}p\rangle \geq (\lambda -L_{i})d^{2}(p,q).
\]
By hypothesis $\lambda > \sup_{i \in I}L_{i}$. Hence, the above equation and Definition \ref{def4.2} imply that $\grad h_{i}$ is  strongly  monotone
with constant $\lambda - \sup_{i \in I}L_{i}$. Therefore, from Proposition \ref{prop:fconv} we conclude that  $h_{i}$ is strongly convex
with constant $\lambda - \sup_{i \in I}L_{i}$.  It easy to see that
\[
\max_{i\in I} h_{i}=f+(\lambda/2) d^2(. \,, \tilde{p} ).
\]
Thus using Proposition \ref{prop:supfcov} the proposition follows.
\end{proof}
\begin{corollary} \label{cor:wdf}
The proximal point method \eqref{E:1.22} applied to $f$ with starting point $p^0=\hat p$ is well defined.
\end{corollary}
\begin{proof}
Assume that $p^k\in L_f(f(q))$ for some $k$. Note that the minimizers of $\psi_k:=f+(\lambda_k/2)d^2(. ,{p^k})$, in case they exist, are in $L_{\psi_k}(\psi_k(p^k))\subset L_f(f(q))\subset\Omega$,  more precisely, 
\[
\argmin_{p\in M}\psi_k(p)=\argmin_{p\in L_{\psi_k}(\psi_k(p^k))}\psi_k(p)).
\]
As $f$ is continuous on $\bar{\Omega}$, $L_{\psi_k}(\psi_k(p^k))$ is closed in the topology of the manifold $M$. Moreover, since $\sup_{\tau\in
T}L_{\tau}<\lambda_k$, we conclude from Lemma~\ref{MPP8} that the application
$\psi_k$ is strongly convex on $\Omega$ with constant $\beta=\lambda-\sup_{\tau\in T}L_\tau$. In this conditions, for $\psi_k$ to has an unique minimizing at $M$ is sufficient that $L_ {\psi_k} (\psi_k (p ^ k)) $ is a bounded set. From the convexity of $\psi_k$ on $\Omega$ we have, in particular, that $L_{\psi_k}(\psi_k(p^k))$ is a convex set. Assume by contradiction that $ L_ {\ psi_k} (\ psi_k (p ^ k)) $ is unbounded. Then, there exists $v\in T_{p^k}M$ and an a geodesic $\gamma:[0,+\infty]\to L_{\psi_k}(\psi_k(p^k))$ such that $\gamma(0)=p^k$ and $\gamma'(0)=v$. Take $t_0>0$. For all $t>t_0$, we have
\[
\psi_k(\gamma(t_0))=\psi_k(\gamma(\frac{t_0}{t}t+(1-\frac{t_0}{t})0)).
\]
Now, taking into account that $\psi_k$ is strongly convex with constant $\beta$, from the last equality follows that 
\[
\psi_k(\gamma(t_0))\leq\frac{t_0}{t}\psi_k(\gamma(t))+(1-\frac{t_0}{t})\psi_k(\gamma(0))-\frac{1}{2}(1-\frac{t_0}{t})\frac{t_0}{t}\beta\|v\|^2t^2.
\]
Because, $\psi_k(\gamma(t)),\psi_k(\gamma(0))\in L_{\psi_k}(\psi_k(p^k)$, from the last inequality, we obtain
\[
\psi_k(\gamma(t_0))\leq \psi_k(p^k)-\frac{1}{2}(t-t_0)t_0\beta\|v\|^2,
\]
which is a contradiction, since $\psi_k(\gamma(t_0))$ is finite and $\psi_k(p^k)-\frac{1}{2}(t-t_0)t_0\beta\|v\|^2$ goes to $-\infty$ as $t$ goes to $+\infty$. Therefore, $L_{\psi_k}(\psi_k(p^k))$ is bounded and, consequently, $p^{k+1}$ is well defined. Since $p^0=\hat p\in L_f(f(q))$, the proof follows from a simple induction argument.
\end{proof}
\begin{lemma}\label{mpprox10}
Let $\{p^k\}$ be the sequence generated by the proximal point method \eqref{E:1.22}. Then the following statements holds:
\begin{itemize}
\item[i)] $f(p^{k+1})+(\lambda_k/2)d^2(p^{k+1},{p^k})\leq f(p^k), \quad k=0, 1, \ldots$;\
\item[ii)] $\{p^{k}\}\subset L_f(f(q))$;
\item [iii)] $0\in\partial \left(f+\frac{\lambda_k}{2}d^2(.\, ,\,{p^k})\right)(p^{k+1}), \quad k=0, 1, \ldots$;
\item[iv)] $-\infty<\bar{f}=\lim_{k \to \infty} f(p^{k})$;
\item[v)] $\lim_{k \to \infty} d(p^{k+1},p^k)=0$.
\end{itemize}
\end{lemma}
\begin{proof}
The first item is an immediate consequence of \eqref{E:1.22}, which implies that $\{f(p^k)\}$ is monotonous and  nonincreasing of where follows the item ii. Since $\max\limits_{i\in I}L_i<\lambda_k$, Lemma~\ref{MPP8} implies $f+(\lambda_k/2)d^2(.\,,p^k)$ convex on $\Omega$, which, together with \eqref{E:1.22}, proof item iii. Using that $\{f(p^k)\}$ is monotonous nonincreasing and that $-\infty<\inf_{p\in M}f(p)$, item iv follows. Finally, item v is a consequence of items i and~iv.
\end{proof}
\begin{lemma} \label{l:conv1}
Let $\{p^k\}$ be the sequence generated by the proximal point method \eqref{E:1.22} with $\lambda_k$ satisfying \eqref{mpp102}. Assume that {\bf h1} and {\bf h2} holds. If $p^k\in L_f(c)$ for some $k$ then $\{p^k\}$ converges to a  point $p^*\in U^*\subset \Omega$.
\end{lemma}
\begin{proof}
By hypotheses, $p^{k}\in L_f(c)$ for some $k$, i.e., there exists $k_0$ such that $f(p^{k_0})\leq c$. Then, from Lemma \ref{mpprox10} item i, $\{p^k\}\subset L_f(c)$ for all $k\geq k_0$. On the other hand, from \eqref{mpp102} we have
\[
\sum_{k=0}^{+\infty}\frac{1}{\lambda_k}=+\infty.
\]
Therefore, using {\bf h1} and {\bf h2} the result follows from similar arguments used in the proof of Theorem 6.1 of \cite{FO2000}.
\end{proof}
\begin{lemma} \label{l:conv2}
Let $\{p^k\}$ be the sequence generated by the proximal point method \eqref{E:1.22} with $\lambda_k$ satisfying \eqref{mpp102}. If  {\bf h3} holds then after a finite number of steps the proximal iterates go into the set $L_f(c)$.
\end{lemma}
\begin{proof}
First note that since $\inf_{p\in M}f(p)<c$ we have $L_f(c)\neq\emptyset$. Suppose by contradiction that $p^k\in L_f(f(q))\setminus L_f(c)$ for all $k$. From \eqref{mpp102} we have $\max\limits_{i\in I}L_i<\lambda_k$. Since  $(1/2)d^2(.,p)$ is a differentiable function with
\[
\grad((1/2)d^2(q,p))=-exp^{-1}_{q}{p},
\]
applying Lemma \ref{MPP8} with $\lambda=\lambda_k$, Lemma \ref{lemacl1} with $g=(1/2)d^2(.,p)$ and the first part of Proposition \ref{SubClarke1}, we obtain:
\[
\partial\left(f+\frac{\lambda_k}{2}d^2(. ,p^k)\right)(p)= \partial^{\circ} f(p)-\lambda_k\exp^{-1}_{p}p^k, \qquad k=0, 1, \ldots.
\]
Hence, it is easily concluded from the last equality and item iii of Lemma \ref{mpprox10} that:
\[
\lambda_k\exp^{-1}_{p^{k+1}}p^k\in \partial^{\circ}f(p^{k+1}), \qquad k=0, 1, \ldots..
\]
Since $p^{k+1}\in L_f(f(q))\setminus L_f(c)$, assumption {\bf h3} and the latter equation give
\[
\|\lambda_k\exp^{-1}_{p^{k+1}}p^k\|>\delta, \qquad k=0, 1, \ldots..
\]
From \eqref{mpp102} we have $\lambda_{k}\leq \bar{\lambda}$. As $d(p^k,p^{k+1})=\|\exp^{-1}_{p^{k+1}}p^k\|$ the last inequality implies that
\[
d(p^k,p^{k+1})>\frac{\delta}{\bar{\lambda}}, \qquad k=0, 1, \ldots.
\]
So, from item v of Lemma \ref{mpprox10} we arrive at a contradiction. Therefore there exists $k_0$ such that $f(p^{k_0})\leq c$ and the result follows from the item i of Lemma \ref{mpprox10}.
\end{proof}
\subsubsection*{Proof of {\bf Theorem \ref{MPP10}}}
\begin{proof}
The well definition of the proximal point method follows from the Corollary \ref{cor:wdf}. Let $\{p^{k}\}$ be the sequence generated by the proximal point method. As $p^0=\hat{p}\in L_f(f(q))$, item~i of Lemma \ref{mpprox10} implies that the whole sequence lies in $L_f(f(q))$.
From item~ iii of Lemma \ref{mpprox10}, we have:
\[
0\in\partial \left(f+\frac{\lambda_k}{2}d^2(.\, ,\,{p^k})\right)(p^{k+1}), \quad k=0, 1, \ldots.
\]
Since $\max\limits_{i\in I}L_i<\lambda_k$, Lemma \ref{MPP8} implies that $f_i+(\lambda_k/2)d^2(.\,,{p^k})$ and $f+(\lambda_k/2)d^2(.\,,{p^k})$ are strongly convex. Thus, applying Corollary~ \ref{CoroMPP1} with  $h_{i}=f_i+(\lambda_k/2)d^2(.\,,{p^k})$ and $h=f+(\lambda_k/2)d^2(.\,,{p^k})$ we conclude that there exists constant $\alpha_i^{k+1}\geq 0$ with $i\in I(p^{k+1})$ such that
\[
0=\sum\limits_{i\in I(p^{k+1})}\alpha_i^{k+1}\grad\left(f_i+\frac{\lambda_k}{2}d^2(.\,,{p^k})\right)(p^{k+1}),  \qquad \sum\limits_{i\in I(p^{k+1})}\alpha^{k+1}_i=1.
\]
This tells us that
\begin{equation}\label{MPP11}
0=\sum\limits_{i\in I(p^{k+1})}\alpha_i^{k+1}\grad f_i(p^{k+1})-\lambda_k\exp^{-1}_{p^{k+1}}p^k,  \qquad \sum\limits_{i\in I(p^{k+1})}\alpha^{k+1}_i=1, \quad k=0, 1, \ldots.
\end{equation}
If the sequence $\{p^k\}$ is finite then there exists $k$ such that $p^{k+1}=p^k$. In this case,   $\exp^{-1}_{p^{k+1}}p^k=0$ and the first  equality in \eqref{MPP11} becomes:
\[
0=\sum\limits_{i\in I(p^{k+1})}\alpha_i^{k+1}\grad f_i(p^{k+1}),
\]
which, together with Proposition \ref{SubClarke1}, implies that $0\in\partial^{\circ}f(p^k)$. Hence, $p^k$ is a stationary point of $f$.

Now, assume that sequence $\{p^k\}$ is infinite and $\bar{p}$ is an accumulation point of it. Let $\{\alpha_i^{k+1}\}\subset \mathbb{R}^{m}$ be the sequence defined by
\[
\alpha^{k+1}=(\alpha_1^{k+1}, \ldots, \alpha_m^{k+1}), \qquad \alpha_j^{k+1}=0, \quad j\in I\backslash I(p^{k+1}).
\]
Since $\sum_{i\in I(p^{k+1})}\alpha^{k+1}_i=1$ we have $\|\alpha^{k+1}\|_{1}=1$ for all $k$, where $\|\;\|_{1}$  denotes the sum norm in $\mathbb{R}^{m}$. Thus $\{\alpha^{k+1}\}$ is bounded. Let $\{p^{k_s+1}\}$ and $\{\alpha^{k_s+1}\}$ be the subsequence of $\{p^{k+1}\}$ and $\{\alpha^{k+1}\}$, respectively, such that $\lim_{s\to +\infty} p^{k_s+1}=\bar{p}$ and $\lim_{s\to +\infty} \alpha^{k_s+1}=\bar{\alpha}$. As $f$ is continuous on $\Omega$, item~ ii of Lemma~ \ref{mpprox10} implies that $\bar{p}\in L_f(f(q))\subset \Omega$. Since $I$ is finite we can assume without loss of generality that
\begin{equation} \label{eq:eii}
I(p^{k_1+1})=I(p^{k_2+1})=...=\bar{I},
\end{equation}
and equation \eqref{MPP11} becomes
\[
0=\sum\limits_{i\in \bar{I}}\alpha_i^{k_{s}+1}\grad f_i(p^{k_{s}+1})-\lambda_{k_{s}}\exp^{-1}_{p^{k_{s}+1}}p^{k_s},  \qquad \sum\limits_{i\in \bar{I}}\alpha^{k_{s}+1}_i=1, \quad s= 1, 2, \ldots.
\]
On the other hand, item~ v of Lemma~ \ref{mpprox10} implies
\[
\lim_{s\to\infty}d(p^{k_{s}+1},p^{k_{s}})=\lim_{s\to\infty}\|\exp^{-1}_{p^{k_{s}+1}}p^{k_{s}}\|=0.
\]
As $\lambda_{k_s}$ is bounded, $\lim_{s\to +\infty} p^{k_s+1}=\bar{p}$ and $\lim_{s\to +\infty} \alpha^{k_s+1}=\bar{\alpha}$, letting $s$ goes to $+\infty$ in the above equality, we conclude
\[
0=\sum_{i\in\bar{I}}\bar{\alpha}_i\grad f_i(\bar{p}), \qquad \sum_{i\in\bar{I}}\bar{\alpha}_i=1.
\]
Using definition of $I(\bar{p})$, equation  \eqref{eq:eii} and the continuity of $f$, we obtain $\bar{I}\subset I(\bar{p})$. Therefore, as $\bar{p}\in \Omega$, it follows from the Proposition~\ref{SubClarke1} that
\[
0\in\partial^{\circ}f(\bar{p}),
\]
i.e., $\bar{p}$ is a stationary point of $f$, which concludes the proof of the first part of the theorem.

The second part follows from Lemma \ref{l:conv1} and Lemma \ref{l:conv2}.
\end{proof}

\section{Examples}\label{sec6}
In this section we present two examples. In the first we consider a non-convex minimization problem where the objective function is defined on a Hadamard manifold with curvature identically zero. In the next example we "generalize" this example to one where the curvature of the Hadamard manifold is not identically zero. In both examples, the classical local proximal point method (see \cite{Kaplan1998}) as well as the Riemannian proximal point method (see \cite{FO2000}) does not apply. However, the method proposed in this paper applies.

\subsection{Example}
Let $(\mathbb{R}_{++}, \langle \, , \, \rangle)$ be the
Riemannian manifold, where
$\mathbb{R}_{++}=\{x\in\mathbb{R}:x>0\}$ and $\langle \, , \,
\rangle$ is the Riemannian metric   $\langle u , v \rangle=g(x)uv$
with $g:\mathbb{R}_{++}\to (0,+\infty)$. The Christoffel
symbol and the geodesic equation are given by
\[
\Gamma(x)=\frac{1}{2}g^{-1}(x)\frac{dg(x)}{dx}=\frac{d}{dx}\ln\sqrt{g(x)}, \qquad \frac{d^2x}{dt^2}+\Gamma(x)\left( \frac{dx}{dt}\right)^2=0,
\]
respectively. Moreover, in relation to the twice differentiable function $h:\mathbb{R}_{++}\to\mathbb{R}$, the Gradient and the Hessian of $h$ are given by
\[
\grad h=g^{-1}h',  \qquad {\rm hess}\ h=h''-\Gamma h',
\]
respectively, where $h'$ and $h''$ denote  the first and second derivatives of $h$ in the Euclidean sense. For more details see \cite{U94}.
In the particular case of $g(x)=x^{-2}$,
\begin{equation}\label{Hess:1}
\Gamma(x)=-x^{-1},\quad \grad h(x)=x^2h'(x),\quad {\rm hess}\ h(x)=h''(x)+x^{-1}h'(x).
\end{equation}
Moreover, the map $\varphi :\mathbb{R} \to \mathbb{R}_{++}$ defined  by $\varphi(x)={\rm e}^x$ is
an isometry between the Euclidean space $\mathbb{R}$ and the manifold $(\mathbb{R}_{++}, \langle \, , \, \rangle)$ and the Riemannian distance $d:\mathbb{R}_{++}\times\mathbb{R_{++}}\to\mathbb{R}_{+}$ is given by
\begin{equation}\label{dRiem:1}
d(x,y)=|\varphi^{-1}(x)-\varphi^{-1}(y)|=|\ln\frac{x}{y}|,
\end{equation}
see, for example \cite{XFLN2006}. Therefore,  $(\mathbb{R}_{++}, \langle \, , \, \rangle)$ is a Hadamard manifold and the unique geodesic $x:\mathbb{R}\to\mathbb{R}_{++}$ with initial conditions $x(0)=x_0$ and $x'(0)=v$ is given by
\[
x(t)=x_0 {\rm e}^{(v/x_0)t}.
\]
From the above equation it is easily seen that any interval  $I\subset \mathbb{R}_{++}$ is a convex set of the manifold $(\mathbb{R}_{++}, \langle \, , \, \rangle)$.

Let $f_1,f_2,f:\mathbb{R}_{++}\to\mathbb{R}$ respectively be given by
\[
f_1(x)=\ln(x),\quad f_2(x)=-\ln(x)+{\rm e}^{-2x}-{\rm e}^{-2},\quad f(x)=\max_{j=1,2}f_j(x),
\]
and consider the problem
\begin{equation}\label{PNCOctPos}
\begin{array}{clc}
   & \min f(x) \\
   & \textnormal{s.t.}\,\,\, x\in\mathbb{R}_{++}.\\
\end{array}
\end{equation}
Take a sequence $\{\lambda_{k}\}$ satisfying $0<\lambda_k$. From \eqref{dRiem:1}, the proximal point method \eqref{E:1.22} becomes
\[
x^{k+1}:=\argmin_{x\in\mathbb{R_{++}}}\left\{f(x)+\frac{\lambda_k}{2}\ln^2\left(\frac{x}{x^k}\right)\right\}, \qquad k=0, 1, \ldots.
\]
Note that $-\infty<\inf_{x\in\mathbb{R}_{++}}f(x)=0$ and, being $f_1$ and $f_2$ twice differentiable functions on $\mathbb{R}_{++}$, the last expression in \eqref{Hess:1} implies that
\begin{equation}\label{ExHess:1}
{\rm hess}\ f_1(x)=0\quad \mbox{and} \quad {\rm hess}\ f_2(x)=(4-\frac{2}{x}){\rm e}^{-2x}, \qquad x\in \mathbb{R}_{++}.
\end{equation}
Let $0<\epsilon< 1/4$, $q=5/16$  and   $\Omega=(\epsilon , +\infty)$. So, $0=\inf_{x\in\mathbb{R}_{++}}f(x)<f(q)$ and $L_f(f(q))\subset\Omega$. Moreover, ${\rm hess}\ f_1$, ${\rm hess}\ f_2$ are bounded on $\Omega$ and consequently $\grad f_1,\grad f_2$ are Lipschitz on $\Omega$. We denote by $L_i$ the constant of Lipschitz of $\grad f_i$, $i=1,2$. Clearly the assumption {\bf h1} of the Theorem \ref{MPP10} is verified with $U^*=\{1\}$.

We claim that there exists $c\in (0,f(1/2))=(0,f_1(1/2))$ such that $L_f(c)$ is convex and $f$ is convex on $L_f(c)$ (in the Riemannian sense). Indeed, as ${\rm hess}\ f_1\geq 0$ in $\mathbb{R}_{++}$ and ${\rm hess}\ f_2\geq 0$ in $[1/2,+\infty)$,  Theorem 6.2 of \cite{XFLN2006} implies that $f_1$ is convex on $\mathbb{R}_{++}$ and $f_2$ is convex on $[1/2,+\infty)$. Thus, it follows from Proposition \ref{prop:supfcov} that $f$ is convex on $[1/2,+\infty)$. Note that for all $c\in (0,f(1/2))$, we have $L_f(c)\cap [1/2,+\infty)=L_f(c)$. Hence, from the convexity of $f$ on $(1/2,+\infty)$ we conclude that $L_f(c)$ is convex, which proves the claim. So, $L_f(c)$ and $f$ satisfy assumption {\bf h2} of Theorem~\ref{MPP10}, for example with $c=f(3/4)$.

Now, note that $L_f(f(q))\backslash L_f(c)=[5/16,3/4)\cup (a, b]$, where $a=(4/3){\rm e}^{({\rm e}^{-3/2}-{\rm e}^{-2})}$ and $b=(16/5){\rm e}^{({\rm e}^{-5/8}-{\rm e}^{-2})}$. Moreover,  $f$ is differentiable on $L_f(f(q))\backslash L_f(c)$ with  $\grad f(x)=\grad f_1(x)$ for $x\in[a, b]$ and $\grad f(x)=\grad f_2(x)$ for $x\in [5/16, 3/4]$. From the second expression in \eqref{Hess:1}, we have:
\[
\grad f(x)=x,\quad x\in (a,b]\quad\mbox{and}\quad \grad f(x)=-x-2x^2{\rm e}^{-2x},\quad x\in [5/16,3/4].
\]
Thus, we have $\|\grad f(x)\|\geq\|\grad f(5/16)\|>2/5$ and $f$ satisfies the assumption {\bf h3} of Theorem \ref{MPP10}.

Summarizing, all assumptions of Theorem 4.1 are satisfied  with $\Omega=(\epsilon , +\infty)$, $q=5/16$, $c=f(3/4)$ and $\delta=2/5$. Therefore, letting $x^0\in \mathbb{R}_{++}$ and $\bar{\lambda}>0 $ such that $x^0\in L_f(f(q))$ and $\max\limits_{i\in I}L_i<\lambda_k\leq\bar{\lambda}$, the proximal point method,  may be applied for solving the above nonconvex problem.

\begin{remark}
Function $f(x)=\max\{\ln(x), -\ln(x)+{\rm e}^{-2x}-{\rm e}^{-2}\}$, in the above example, is  nonconvex  (in the Euclidean sense) when restricted to any open neighborhood containing its minimizer $x^*=1$. Therefore, the classical local proximal point method (see \cite{Kaplan1998}) cannot be applied to minimize this function. Also, as $f$ is nonconvex in the Riemannian sense as well, the Riemannian proximal point method (see \cite{FO2000}) cannot be applied to minimize this function either.
\end{remark}

\subsection{Example}\label{sec7} Let ${\mathbb S}^{n}$ be the set
of the symmetric matrices, ${\mathbb S}^{n}_{+}$ be the cone of
the symmetric positive semi-definite matrices and ${\mathbb
S}^{n}_{++}$ be the cone of the symmetric positive definite
matrices both $n\times n$. For $X, \, Y\in {\mathbb S}^{n}_{+}$,
$Y\succeq X$ (or $X \preceq Y$) means that $Y-X \in {\mathbb
S}^{m}_{+}$ and $Y\succ X$ (or $X \prec Y$) means that $Y-X \in
{\mathbb S}^{n}_{++}$. We will denote the Frobenius norm by $\|\,.\,\|_F$.

Following Rothaus \cite{Rot1960}, let $M:=({\mathbb S}^n_{++}, \langle \, , \, \rangle)$ be the Riemannian manifold endowed with the Riemannian metric induced by the Euclidean Hessian of $\Psi(X)=-\ln\det X$,
\begin{equation}\label{Rdt:0.0}
\langle U,V \rangle=\tr (V\Psi''(X)U)=\tr (VX^{-1}UX^{-1}),\qquad X\in M, \qquad U,V\in
T_XM,
\end{equation}
where $\tr(A)$ denotes the trace of matrix $A\in {\mathbb S}^n$ and $T_XM\approx~\mathbb{S}^n$, with the corresponding norm denoted by ~ $\|\; .\; \|$. In this case the unique geodesic segment connecting any $X,Y\in M$ is given by
\[
\gamma(t)=X^{1/2}\left(X^{-1/2}YX^{-1/2}\right)^tX^{1/2},\qquad t\in[0,1],
\]
see, for instance, \cite{NT01}.	
More precisely, $M$ is a Hadamard manifold, see for example \cite{Lang 1998}, Theorem 1.2. page 325. From above equality it is immediate that
\[
\gamma'(0)=X^{1/2}\ln \left(X^{-1/2}YX^{-1/2}\right)X^{1/2}.
\]
Thus, for each $X\in M$, $\exp^{-1}_X:M\to T_XM$ and $\exp_X:T_XM\to M$ are given, respectively, by
\begin{equation}\label{Exp:1}
\exp^{-1}_XY=X^{1/2}\ln \left(X^{-1/2}YX^{-1/2}\right)X^{1/2},\qquad \exp_XV=X^{1/2}e^{\left(X^{-1/2}YX^{-1/2}\right)}X^{1/2}.
\end{equation}
Now, since the Riemannain distance $d$ is given by $d({X}\, , \, Y)\,=\,||exp^{-1}_{X}Y||$, from \eqref{Rdt:0.0} along with first expression in \eqref{Exp:1}, we conclude that:
\begin{equation}\label{dRiem:1000}
d^2(X,Y)=\tr\left(\ln^2 X^{-1/2}YX^{-1/2}\right)=\sum^n_{i=1}\ln^2\lambda_i\left(X^{-\frac{1}{2}}YX^{-\frac{1}{2}}\right),
\end{equation}
where $\lambda_i(X^{-\frac{1}{2}}YX^{-\frac{1}{2}})$ denotes the $i^{th}$ eigenvalue of the matrix
$X^{-\frac{1}{2}}YX^{-\frac{1}{2}}$. The gradient and the Hessian of a twice differentiable function $F:{\mathbb S}^n_{++} \to \mathbb{R}$ is given, respectively,  by:
\[
\mbox{grad} F(X)=XF'(X)X, \
\mbox{hess}\,F(X)\big{(}V, V\big{)}=\mbox{tr}\big{(} VF'' (X) V\big{)}+
\mbox{tr}\big{(}F'(X)VX^{-1}V \big{)},\   V\in T_XM,
\]
where $F'(X)$ and  $F''(X)$ are the  Euclidean gradient and Hessian, respectively. We remind that
a twice differentiable function $F:{\mathbb S}^n_{++}\to \mathbb{R}$ is convex on the manifold $M$ if it satisfies the
condition:
\[
\mbox{hess}\,F(X)\big{(}V, V\big{)}\geq 0,\qquad \qquad X\in M, \quad V\in T_XM.
\]
Let $F_1,F_2,F_3:{\mathbb S}^n_{++}\to\mathbb{R}$ be given, respectively, by:
\[
F_1(X)=\ln\det X ,\quad F_2(X)=-4\ln\det X+e^{-2\tr X}-e^{-2n},\quad F_3(X)=\tr X^{-1}-n.
\]
Note that $F_1$, $F_2$ and $F_3$ are twice differentiable functions on ${\mathbb S}^n_{++}$. Taking $X\in T_XM$ and $V\in T_XM$, the Euclidean Gradients are given by:
\[
F'_1(X)=X^{-1},\quad F'_2(X)=-4X^{-1}-2e^{-2\tr X},\quad F'_3(X)=-X^{-2},
\]
and the Euclidean Hessians are given by:
\begin{align*}
F''_1(X)V&=-X^{-1}VX^{-1}, \\
F''_2(X)V&=4X^{-1}VX^{-1}+4e^{-2\tr X}V,\\
F''_3(X)V&=X^{-1}VX^{-2}+X^{-2}VX^{-1}.
\end{align*}
Thus, the Riemannian gradients of $F_1, F_2$ and $F_3$ are given by:
\begin{equation}\label{Rdt:0}
\grad F_1(X)=X,\quad \grad F_2(X)=-4X-2e^{-2\tr X}X^2,\quad \grad F_3(X)=-I,
\end{equation}
where $I$ denotes the identity  matrix, and Riemannian Hessians of $F_1, F_2$ and $F_3$ are given by:
\begin{align}\label{GLip:1}
\mbox{hess}\,F_{1}(X)\left(V, V\right)&=0,\nonumber \\
\mbox{hess}\,F_{2}(X)\left(V, V\right)&=2e^{-2\tr X}\tr\left[(2I-X^{-1})V^2\right],\\
\mbox{hess}\,F_{3}(X)\left(V, V\right)&=\tr\left(X^{-1}VX^{-2}V\right)=\|X^{-1}VX^{-1/2}\|_F^2\nonumber.
\end{align}
It is easily seen that the functions $F_{1}$ and $F_{3}$ are convex everywhere and the function $F_{2}$ is convex on any convex subset of the set
\begin{equation} \label{eq:c}
C:=\left\{X\in {\mathbb S}^n_{++}:2I-X^{-1}\succ 0\right\}=\left\{X\in {\mathbb S}^n_{++}:\lambda_{\min}(X)> 1/2\right\},
\end{equation}
where $\lambda_{\min}(A)$ denotes the minimum eigenvalue of  the matrix $A$.

Let $F:{\mathbb S}^n_{++}\to\mathbb{R}$ be given by
\[
F(X)=\max_{j=1,2,3}F_j(X).
\]
From Proposition \ref{prop:supfcov} $F$ is convex on any convex subset of $C$. Consider  the following optimization problem
\begin{equation}\label{PNCMatrix}
\begin{array}{clc}
   & \min F(X) \\
   & \textnormal{s.t.}\,\,\, X\in {\mathbb S}^n_{++}.\\
\end{array}
\end{equation}
If $\{\eta_{k}\}$ is a sequence satisfying $0<\eta_k$, then from \eqref{dRiem:1000} the proximal point method \eqref{E:1.22} becomes
\[
X_{k+1}={\arg\min}_{X \in {\mathbb S}^n_{++}}\left \{ f(X)+
\frac{\eta_k}{2}\sum^n_{i=1}\ln^2
\lambda_i\left(X^{-\frac{1}{2}}X_{k}X^{-\frac{1}{2}}\right)\right \}.
\]
Define the sets $U_1, U_2$ and $U_3$ as
\[
U_1:=\{X\in {\mathbb S}^n_{++}:\det X=1\},\ U_2:=\{X\in {\mathbb S}^n_{++}:\det X>1\},\ U_3:=\{X\in {\mathbb S}^n_{++}:\det X<1\},
\]
and note that ${\mathbb S}^n_{++}=U_1 \cup U_2 \cup U_3$ and $U_{i}\cap U_{j}=\emptyset,$ for $i, j=1,2,3$, $i\neq j$.
\begin{claim}\label{claim6:1}
\begin{equation}\label{eknp:1}
F(x)=\begin{cases}
F_3(X)\geq 0, &\; X\in U_1;\\
\max\{ F_1(X),F_3(X)\}>0,  &\; X\in U_2;\\
\max\{ F_2(X),F_3(X)\}>0,  &\; X\in U_3.
\end{cases}
\end{equation}
In particular, $\inf_{{\mathbb S}^n_{++}}F(X)=0$ and $U^*=\{I\}$, where $U^*$ is the solution set of the problem~ \eqref{PNCMatrix} and $I$ is the $n\times n$ identity matrix.
\end{claim}
\begin{proof}
First of all, note that for each $X\in {\mathbb S}^{n}_{++}$, we have
\begin{equation}\label{Rdt:1}
\sqrt[n]{\det X}\leq\frac{\tr X}{n}.
\end{equation}
If  $X\in U_1$ then  $\det X=1$. Using the above inequality we obtain  $\tr X\geq n$ which implies
\[
F_2(X)=e^{-2\tr X}-e^{-2n}\leq 0=F_1(X).
\]
On the other hand $\det X^{-1}=1$. Hence, \eqref{Rdt:1} gives us $\tr X^{-1}-n\geq 0$ and  thus $F_3(X)\geq 0$. So, from definition of $F$ and taking into account the above inequality, we conclude
\begin{equation}\label{Rdt:2}
F(X)=F_3(X)\geq 0,\qquad  X\in U_1.
\end{equation}
If $X\in U_2$ then $\det X>1$. Thus,  from  \eqref{Rdt:1} we also obtain  $\tr X>n$. Hence,
\[
F_2(X)=-4\ln\det X+e^{-2\tr X}-e^{-2n}<-\ln\det X< 0.
\]
Since, $\det X>1$ we have $F_1(X)=\ln\det X>0$. Therefore,
\begin{equation}\label{Rdt:3}
F(X)=\max\{F_1(X),F_3(X)\}>0, \qquad \;X\in U_2.
\end{equation}
Finally, if $X\in U_3$ then   $\det X^{-1}>1$. Hence, inequality  \eqref{Rdt:1} implies that $\tr X^{-1} > n$ and consequently $F_3(X)=\tr X^{-1}-n>0$. As $\det X<1$ we have $F_1(X)=\ln\det(X)<~0$. Thus,
\begin{equation}\label{Rdt:4}
F(X)=\max\{F_2(X),F_3(X)\}>0,\qquad \; X\in U_3,
\end{equation}
and the the first part of the claim is proved.

From the first part we conclude that $F(X)\geq 0$, for all $X\in {\mathbb S}^n_{++}$. Since $F(I)=0$  we have $\inf_{{\mathbb S}^n_{++}}F(X)=0$. In order to prove the last statement,  note that  if $F(X)=0$ then $X\in U_1$. Definition of $U_1$  and  \eqref{eknp:1} give us  $\det X^{-1}=1$ and  $F(X)=\tr(X^{-1})-n$  for all  $X\in U_1$. Hence, if $X\in U^*$ then
$\tr(X^{-1})=n$ and $\det X^{-1}=1$, which implies that \eqref{Rdt:1} holds with equality. On the other hand,  \eqref{Rdt:1} holds with equality if only if $\lambda_1(X)=...=\lambda_n(X)>0$. As $\tr(X^{-1})=n$ we conclude that $\lambda_i(X)=1$ for $i=1, \ldots, n$ and the result follows.
\end{proof}
\begin{claim}\label{claim6:2}
$F_1(X)> F_3(X)$ for $X\in A$, where
\[
A:=\left\{X\in {\mathbb S}^n_{++}:X-I\succ 0\right\}=\left\{X\in {\mathbb S}^n_{++}:\lambda_{\min}(X)>1\right\}\subset U_2.
\]
Consequently, $F(X)=F_1(X)$ for $X\in A$.
\end{claim}
\begin{proof}
Define the following   function
\[
\psi(X):=F_1(X)-F_3(X).
\]
From the Mean Value Theorem on the Euclidean space ${\mathbb S}^{n}_{++}$ there exists $\tilde{X}\in U_2$ such that
\[
\psi(X)=\psi(I)+\tr\left(\psi'(\tilde{X})(X-I)\right)=\tr\left((\tilde{X}^{-1}+\tilde{X}^{-2})(X-I)\right).
\]
If $X\in A$ then $X-I\in {\mathbb S}^{n}_{++}$. Since $X-I$ and  $\tilde{X}^{-1}+\tilde{X}^{-2}$ belong  to ${\mathbb S}^{n}_{++}$ and the trace of the product of positive definite matrices is positive, we conclude from the above equality that $\psi(X)=F_1(X)-F_3(X)>0$ for all $X\in A$, which proof the first part of the claim. The second part of the claim follows  from Claim~\ref{claim6:1} .
\end{proof}
\begin{claim}\label{claim6:3}
$F(X)=F_2(X)$ for $X\in B$, where
\[
B:=\left\{X\in {\mathbb S}^n_{++}:(1/4)I\prec X\prec I\right\}=\left\{X\in {\mathbb S}^n_{++}:1/4<\lambda_{\min}(X),\, \lambda_{\max}(X)<1\right\}.
\]
\end{claim}
\begin{proof}
Note that $B\subset U_3\cap\{X\in {\mathbb S}^n_{++}:\tr X \leq n\}$. So, from \eqref{Rdt:4}, to prove this claim it suffices to verify that $F_2(X)>F_3(X)$ for all $X\in B$. Let $\phi:\mathbb{R}_{++}\to\mathbb{R}$ be the function given $\phi(t):=-4\ln t -1/t +1$. From the convexity of $-4\ln t$, we have
\begin{equation} \label{Rdt:4.5}
\phi(t)\geq -4\ln 1-4(t-1)-\frac{1}{t}+1> 0, \qquad t\in (1/4 \,,1).
\end{equation}
On the other hand,  definitions of $F_2$, $F_4$ and $\phi$ yields:
\[
F_2(X)-F_4(X)=e^{-2\tr X}-e^{-2n}+\sum^{n}_{j=1}\phi(\lambda_j(X)).
\]
Since $e^{-2\tr X}-e^{-2n}> 0$ for $X\in B$, combining \eqref{Rdt:4.5} with the latter equality  we conclude  that  $F_2(X)>F_3(X)$ for $X\in B$. Hence, the desired equality follows.
\end{proof}
It is easy to verify that $F$ is a coercive function. Thus, $L_F(a)$ is bounded for all $a\in\mathbb{R}$, i.e., all the level sets of $F$ are bounded. Take
\begin{equation} \label{eq:q}
Q:=\diag(1/4, \ldots, 1/4).
\end{equation}
As $L_F(F(Q))\subset {\mathbb S}^n_{++}$ is a bounded set, then there exists $\kappa>4$ sufficiently large so that
\begin{equation}\label{Rdt:5}
\lambda_{\min}(X)>1 / {\kappa},\quad \lambda_{\max}(X)<{\kappa},\qquad\; X\in L_F(F(Q)).
\end{equation}
Now we  define
\[
\Omega:=\{X\in {\mathbb S}^n_{++}:\tr X^{-1}< \kappa n\}.
\]
Since the function ${\mathbb S}^n_{++}\ni X\mapsto \tr X^{-1}$ is  convex, the set  $\Omega$ is convex. From Claim ~ \ref{claim6:1} $0=\inf_{{\mathbb S}^n_{++}}F(X)<F(Q)$ and, consequently, $\inter L_F(F(Q))\neq\emptyset$, where $\inter A$ represents the interior of the set $A$. Moreover, it is immediate to verify that
\[
C\subset\Omega, \qquad I\in \inter L_F(F(Q))\subset L_F(F(Q))\subset \Omega.
\]
\begin{claim} \label{claim:L}
The gradient  vector fields $\grad F_1, \grad F_2$ and $\grad F_3$ are Lipschitz on $\Omega$.
\end{claim}
\begin{proof}
From the first equality in \eqref{GLip:1} and Proposition \ref{gradLips:01} is immediate
that the gradient vector field $\grad F_1$ is Lipschitz. Now, the second equality in \eqref{GLip:1} implies that:
\[
|\mbox{hess}\,F_{2}(X)\left(V, V\right)|=2e^{-2\tr X}|\tr\left[(2I-X^{-1})V^2\right]|,\qquad X\in M,\; V\in T_XM.
\]
As $\tr X>0$ and using the Cauchy-Schwartz inequality, we obtain from the above equality that:
\[
|\mbox{hess}\,F_{2}(X)\left(V, V\right)|\leq 2\|2I-X^{-1}\|_F\|V^2\|_F,\qquad X\in M,\; V\in T_XM
\]
Hence, using the definition of the metric we conclude that:
\[
|\mbox{hess}\,F_{2}(X)\left(V, V\right)|\leq 2\left(2\sqrt{n}+\|X^{-1}\|_F\right)\|V^2\|_F,\qquad X\in M,\; V\in T_XM.
\]
Thus, as $\|X^{-1}\|_F\leq \tr X^{-1}$, from the latter inequality along with the definition of $\Omega$, we conclude that:
\begin{equation}\label{gradLipsch:1}
|\mbox{hess}\,F_{2}(X)\left(V, V\right)|\leq 4\sqrt{n}+2\kappa n,\qquad X\in\Omega,\; \|V\|_F=1.
\end{equation}
On the other hand, using the third equality in \eqref{GLip:1} and the Cauchy-Schwarz inequality, we obtain:
\[
|\mbox{hess}\,F_{3}(X)\left(V, V\right)|=|\tr \left(X^{-1}VX^{-2}V\right)|\leq \|X^{-1}VX^{-1}\|_F\|VX^{-1}\|_F,\quad X\in M,\; V\in T_XM.
\]
Since $\|.\|_F$ is submultiplicative and reusing that $\|X^{-1}\|_F\leq \tr X^{-1}$ above inequality implies that:
\[
|\mbox{hess}\,F_{3}(X)\left(V, V\right)|\leq \left(\tr X^{-1}\right)^3\|V\|_F^2,\quad X\in M,\; V\in T_XM.
\]
Thus,
\begin{equation}\label{gradLipsch:2}
|\mbox{hess}\,F_{3}(X)\left(V, V\right)|\leq \kappa^3n^3,\qquad X\in\Omega,\; \|V\|_F=1.
\end{equation}
From \eqref{gradLipsch:1} and \eqref{gradLipsch:2} we conclude that $hess\,F_{2}$ and $hess\,F_{3}$ are also bounded operators when restricted to $\Omega$. Therefore, from Proposition \ref{gradLips:01}, $\grad F_1$, $\grad F_2$ and $\grad F_3$ are Lipschitz on $\Omega$.
\end{proof}
We denote by $L_i$ the Lipschitz constant  of $\grad f_i$, $i=1,2,3$ on $\Omega$.

From Claim \ref{claim6:1} the assumption {\bf h1} of the Theorem \ref{MPP10} is satisfied  with $U^*=\{I\}$. The next result will be useful  to assure that the assumption {\bf h2} of the Theorem \ref{MPP10} is satisfied.

\begin{claim}\label{claim6:4}
There exists $\sigma>0$ and $\hat{X}\in {\mathbb S}^n_{++}\backslash\{I\}$ such that,
\[
L_F(F(\hat{X}))\subset B_{\sigma}(I)\subset \inter L_F(F(Q))\cap  C,
\]
where $B_{\sigma}(I):=\{X\in {\mathbb S}^n_{++}: d(X,I)<\sigma\}$, $C$ is  defined in  \eqref{eq:c} and $Q$ in \eqref{eq:q}.
\end{claim}
\begin{proof}
Since $U^*=\{I\}$ and $I\neq Q$, we conclude that $I\in\inter L_F(F(Q))\cap C$. So, it is easily seen that there exists $\sigma>0$ such that:
\[
B_{\sigma}(I)\subset\inter L_F(F(Q))\cap C,
\]
proving the last inclusion. In order to prove the first inclusion, assume by contradiction that for every $X\in B_{\sigma}(I)$ we have
$
L_F(F(X))\cap (\mathbb{S}^{n}_{++}\backslash B_{\sigma}(I))\neq\emptyset.
$
Therefore, there exist  sequences  $\{X^k\}, \, \{Y^k\}\subset \mathbb{S}^n_{++}$ such that
\begin{equation} \label{Rdt:5.1}
X^k\in  B_{\sigma}(I), \quad \lim_{k\to \infty}X^k=I, \qquad Y^k  \in L_F(F(X^k))\backslash B_{\sigma}(I), \qquad k=0, 1, \ldots.
\end{equation}
Since $\{X^k\}\in B_{\sigma}(I)\subset  L_F(F(Q))$ we have  $L_F(F(X^k))\subset L_F(F(Q))$ for $k=0, 1, \ldots$. Hence, $\{Y^k\} \subset L_F(F(Q))$. As $L_F(F(Q))$ is a bounded set, we  assume  (taking a subsequence, if necessary) that $\{Y^k\}$ converges to some $\bar{Y}\in L_F(F(Q))$. Thus, continuity of $F$ along with the two last equations in  \eqref{Rdt:5.1} imply that:
\[
F(\bar{Y})\leq F(I), \quad  \quad \bar{Y}\neq I,
\]
and, since $U^*=\{I\}$, we obtain a contradiction, which proves the claim.
\end{proof}
From the above claim along with the convexity of the ball $B_\sigma(I)$, we conclude that the assumption {\bf h2} is satisfied  with $c=F(\hat X)$.

\begin{claim}\label{claim6:5}
There exist $\tilde{\delta}>0$ such that
\begin{equation}\label{Rdt:6}
\|\grad F_j(X)\|>\tilde{\delta}, \qquad \forall\; X\in L_F(F(Q)), \qquad j=1, 2, 3,
\end{equation}
where $Q$ is defined in \eqref{eq:q}.
\end{claim}
\begin{proof}
Let $X\in M$. Using \eqref{Rdt:0} along with the definition of the Riemannian metric, we have:
\begin{align*}
\|\grad F_1(X)\|&=\sqrt{n},\\
\|\grad F_2(X)\|&=\sqrt{\tr(I+4e^{-2\tr X}\tr X+4e^{-4\tr X}\tr X^2)},\\
\|\grad F_3(X)\|&=\sqrt{\tr X^{-2}}.
\end{align*}
Now, if $X\in L_F(F(Q))$, then from the second inequality in \eqref{Rdt:5} we conclude that:
\[
\|\grad F_3(X)\|>\sqrt{n}/\kappa.
\]
It is easily seen that $\|\grad F_2(X)\|>\sqrt{n}$. Thus, the result follows by taking $\tilde{\delta}=\sqrt{n}/\kappa$.
\end{proof}
For  checking that $F$ satisfies the assumption {\bf h3} we will use the following notation
\begin{equation} \label{eq:eb}
B^E_r(I):=\{X\in S^n_{++}:\|X-I\|_F<r \},
\end{equation}
for the Euclidean ball with center in $I$ and radius $r$. Since $I\in \inter L_F(F(\hat X))$, where $\hat X$ is as defined in Claim \ref{claim6:4}, take $r>0$ such that
\[
B^E_r(I) \subset L_F(F(\hat X)).
\]
Using Claim \ref{claim6:4} we have  $L_F(F(\hat{X}))\subset L_F(F(Q))$, with $Q$ as in \eqref{eq:q}. Then in order to prove that $F$ satisfies the assumption {\bf h3} it is sufficient to prove:
\begin{claim}
There exists $\delta>0$ such that
\begin{equation}\label{Rdt:6.1}
\|y(X)\|>\delta, \qquad   X\in  L_F(F(Q))\backslash B^E_r(I), \quad  y(X)\in \partial^{\circ}F(X).
\end{equation}
\end{claim}
\begin{proof}
Let  $B^E_r(I)$ be  the Euclidean ball as defined in \eqref{eq:eb}. Take $X\in  L_F(F(Q))\backslash B^E_r(I)$ and  consider the active indexes set
\[
I(X):=\{i : F(X)=F_i(X), i=1,2,3\},
\]
of $F$ at $X$. Hence  Claim \ref{claim6:1} implies that
\begin{equation} \label{ia}
I(X)=\begin{cases}
\{3\}, &\; X\in U_1;\\
\{1\},\;\{3\},\; or\; \{1,3\},  &\; X\in U_2;\\
\{2\},\;\{3\},\; or\; \{2,3\}, &\; X\in U_3.
\end{cases}
\end{equation}
So, we conclude from Claim \ref{claim6:5} that in order to prove the claim, it is suffices to consider the following cases:
\begin{itemize}
\item[a)] $X\in D_1=\{X\in U_2:I(X)=\{1,3\}\}$;
\item[b)] $X\in D_2=\{X\in U_3:I(X)=\{2,3\}\}$.
\end{itemize}
Let us suppose, initially, that $a)$ holds, i.e., $X\in D_1$. Take $y(X)\in \partial^{\circ}F(X)$. Lemma \ref{ap:1}  guarantees that it exists $\alpha\in[0,1]$ such that:
\[
y(X)=\alpha\grad F_1(X)+(1-\alpha)\grad F_3(X).
\]
Now, using \eqref{Rdt:0}, the definition of the Riemannian metric and the latter expression, we have:
\[
\|y(X)\|^2=\tr\left(\alpha I-(1-\alpha)X^{-1}\right)^2,
\]
which after simple algebraic manipulations becomes
\begin{equation}\label{Rdt:7}
\|y(X)\|^2=\sum^n_{j=1}\left( \alpha- (1-\alpha)\lambda^{-1}_j(X)\right)^2,
\end{equation}
where $\lambda_j(X)$ denotes the $j^{th}$ eigenvalue of the matrix $X$. From \eqref{Rdt:7}  it is easy to see that
\[
\|y(X)\|^2\geq \left( \alpha- (1-\alpha)\lambda^{-1}_{\min}(X)\right)^2+\left( \alpha- (1-\alpha)\lambda^{-1}_{\max}(X)\right)^2.
\]
Minimizing the second degree polynomial of the right hand side in the above inequality with respect to $\alpha$, we obtain:
\begin{equation}\label{Rdt:8}
\|y(X)\|^2\geq\frac{\left(\lambda^{-1}_{\min}(X)-\lambda^{-1}_{\max}(X)\right)^2}{\left(1+\lambda^{-1}_{\min}(X)\right)^2+\left(1+\lambda^{-1}_{\max}(X)\right)^2}.
\end{equation}
Since $X\in D_1\subset U_2$, combining Claim~\ref{claim6:1} with Claim~\ref{claim6:2}, we have $\lambda_{\min}(X)\leq 1$. Furthermore, in this case $\det X > 1$ which, combining with $\lambda_{\min}(X)\leq 1$, implies that
\[
\lambda_{\max}(X)>1.
\]
On the other hand, as $\det X > 1$ equation \eqref{Rdt:1}  yields $\tr X>n.$ Since $X\notin B^{E}_r(I)$ and $\tr X>n$, we conclude from above inequality and after some algebraic manipulations that
\[
r^2\leq \tr(X-I)^2=\tr X^2-2\tr X+n\leq n \lambda_{\max}^2(X)-2n+ n\lambda_{\max}^2(X),
\]
which implies that
\[
\lambda_{\max}(X)>\sqrt{1+r^2/2n}.
\]
Combining first the inequality in \eqref{Rdt:5} with the inequalities \eqref{Rdt:8} and the latter inequality, we obtain:
\begin{equation} \label{fp}
\|y(X)\|>\delta_{1}=\frac{1-1/\sqrt{1+r/2n}}{\sqrt{(1+\kappa)^2+(1+1/\sqrt{1+r/2n})^2}}>0, \quad X\in D_1, \quad y(X)\in \partial^{\circ}F(X).
\end{equation}
Now, let us suppose that  the item $b)$ holds, i.e., $X\in D_1$. Take $y(X)\in \partial^{\circ}F(X)$. Due to Lemma \ref{ap:1} there exists $\alpha\in[0,1]$ such that
\[
y(X)=\alpha\grad F_1(X)+(1-\alpha)\grad F_3(X).
\]
Thus, \eqref{Rdt:0}, definition of the Riemannian metric and the last expression imply that
\[
\|y(X)\|^2=\tr\left(\alpha (-X-2e^{-2\tr X}X^2)-(1-\alpha)I\right)^2.
\]
Some algebraic manipulations in the last equality and taking into account that $\lambda_j(X)>0$, for all $j=1,\ldots,n$ yields:
\begin{equation} \label{Rdt8.1}
\|y(X)\|^2>\alpha^2n+2\alpha(1-\alpha)\tr X^{-1}+(1-\alpha)\tr X^{-2}.
\end{equation}
Due to $X\in D_2\subset U_3$ we have $\det X^{-1}>0$.  Hence, using \eqref{Rdt:1}, we obtain:
\[
\tr X^{-1}>n, \qquad \tr X^{-2}>n
\]
Since $\alpha \in [0,1]$, substituting two above inequalities into \eqref{Rdt8.1}  we obtain
\begin{equation} \label{sp}
\|y(X)\|^2>n, \qquad X\in D_2, \quad y(X)\in \partial^{\circ}F(X).
\end{equation}
Therefore, taking  $\delta=\min\{\tilde{\delta},\delta_1,n\}>0$, where $\tilde{\delta}$ is given in the last claim, the claim follows from \eqref{ia}, \eqref{fp}, \eqref{sp} and  Claim~\ref{claim6:5}.
\end{proof}
Thus $L_F(F(Q))\backslash L_F(F(\hat{X}))\subset L_F(F(Q))\backslash B^E_r(I)$, the assumption {\bf h3} is satisfied with $\delta=\min\{\tilde{\delta},\delta_1,\delta_2\}>0$.

Summarizing, all assumptions of Theorem 4.1 are satisfied by taking $\Omega=\{X\in {\mathbb S}^n_{++}:\tr X^{-1}< \kappa n\}$, $q=\diag(1/4, \ldots, 1/4)$, $c=F(\hat X)$, where $\hat X$ is as defined in Claim \ref{claim6:4}, and $\delta>0$ is defined above. Letting $X^0\in \mathbb{S}^n_{++}$ and $\bar{\lambda}>0 $ such that $X^0\in L_f(f(q))$ and $\max\limits_{i\in I}L_i<\lambda_k\leq\bar{\lambda}$, the proximal point method, may be applied for solving the above nonconvex problem \eqref{PNCMatrix}.

\begin{remark}
In agreement with the claims \ref{claim6:1} and \ref{claim6:2}, the function
\[
F(X)=\max\{\ln\det X,\\-4\ln\det X+e^{-2\tr X}-e^{-2n},\tr X^{-1}-n\},
\]
in the above example, is nonconvex (in Euclidean sense) when restricted to any open convex neighborhood containing its minimizer $X^*=I$. Therefore, the classical local proximal point method cannot be applied to minimize this function. From the second relationship in \eqref{GLip:1} together Claim \ref{claim6:3} we conclude that the function $F$ is nonconvex in the Riemannian sense as well. Therefore the Riemannian proximal point method cannot also be applied to minimize this function either.
\end{remark}
\section{Final Remarks}
We have extended the range of application of the proximal point method to solve nonconvex optimization problems on Hadamard manifolds, namely, in the case of the objective function being given by the maximum of a certain class of continuously differentiable functions. We have certified, through examples, that the class of minimization problems for which the local proximal point method can be applied is different from the class of minimization problem for which the  classical local proximal point as well as the Riemannian proximal point method are applied. An interesting subject now is to extend the proximal point method to minimization problems where objective functions are Lower-$C^2$ type.

\section{Appendix}
In this section we will use notation as in Example \ref{sec7}.

\begin{lemma}\label{ap:0.1}
If $X\in M$ and $V\in T_XM$, then
\[
\lim_{t\to 0\ Y\to X}\frac{d\left(\exp_Y t(D\exp_X)_{\exp^{-1}_XY}V,
\,Y+tV\right)}{t}=0.
\]
\end{lemma}
\begin{proof}
Preliminarily, note that
\[
d\left(\exp_Y t(D\exp_X)_{\exp^{-1}_XY}V,
\,Y+tV\right)=d\left(\exp_Yt(D\exp_X)_{\exp^{-1}_XY}V,\exp_Y\left(\exp^{-1}_Y(Y+tV)\right)\right).
\]
Since the exponential map $exp:TM\to M$ is smooth, it is in particular locally Lipschitz. Let $U_X\subset M$ be a neighborhood of $X$ such that $TU_X\approx U_X\times T_XM$ and the map $exp$ is Lipschitz
in $TU_X$ with constant $K$. Take $\delta=\delta(X)>0$ such that, for all $Y\in B_{\delta}(X)$ and $t\in(0,\delta)$,
\begin{equation}\label{eq:linq}
d\left(\exp_Yt(D\exp_X)_{\exp^{-1}_XY}V,\exp_Y\left(\exp^{-1}_Y(Y+tV)\right)\right)\leq K\|{\cal H}(t,Y)\|,
\end{equation}
where ${\cal H}(t,Y)=tD\exp_X)_{\exp^{-1}_XY}V-\exp^{-1}_Y(Y+tV)$.
On the other hand, by the definition of the metric \eqref{Rdt:0.0}
\[
\|{\cal H}(t,Y)\|^2=\tr\left(\left(t(D\exp_X)_ {\exp^{-1}_XY}V-\exp^{-1}_Y(Y+tV) \right)Y^{-1}\right)^2.
\]
The definition of $\exp^{-1}_{Y}$ in \eqref{Exp:1} and above equality yields:
\[
\|{\cal H}(t,Y)\|^2\leq\tr\left[\left(t(D\exp_X)_ {\exp^{-1}_XY}V-Y^{1/2}\ln\left(I+tY^{-1/2}VY^{-1/2}\right)Y^{1/2}\right)Y^{-1}\right]^2,
\]
which, after simple algebraic manipulation, implies:
\begin{equation}\label{drcms:2}
\left\|\frac{{\cal H}(t,Y)}{t}\right\|^2\leq\tr\left[\left((D\exp_X)_ {\exp^{-1}_XY}V-Y^{1/2}\ln\left(I+tY^{-1/2}VY^{-1/2}\right)^{1/t}Y^{1/2}\right)Y^{-1}\right]^2.
\end{equation}
It is easily seen, that:
\[
\lim_{Y\to X}D\exp_X)_ {\exp^{-1}_XY}V=V, \qquad \lim_{t\to 0 \, Y\to X}\ln\left(I+tY^{-1/2}VY^{-1/2}\right)^{1/t}=X^{-1/2}VX^{-1/2}.
\]
Hence, combining inequality \eqref{drcms:2} with the two latter equalities, we conclude that:
\[
\lim_{t\to 0 \, Y\to X}\left\|\frac{{\cal H}(t,Y)}{t}\right\|=0.
\]
The, last equality along with \eqref{eq:linq} imply the desired result.
\end{proof}
\begin{lemma}\label{ap:0}
Let $\Omega\subset \mathbb{S}^n_{++}$ be an open convex set. If  $F$ is a locally  Lipschitz function on  $\Omega$, then
\[
F^{\circ}(X,V)=F_E^{\circ}(X,V), \qquad X\in \Omega, \qquad  V\in T_XM.
\]
\end{lemma}
\begin{proof}
Take $X\in \Omega$ and $V\in T_XM$. Since $F$ is locally  Lipschitz on $\Omega$ and $TM$ is locally a product, there exists $\delta=\delta(X)>0$ such that $TB_{\delta}(X)\approx B_{\delta}(X)\times\mathbb{R}^n$ and
\[
\left|\frac{F({\cal G}(t, Y))-F(Y+tV)}{t}\right| \\ \leq L_X
\,\frac{d({\cal G}(t, Y), \,Y+tV)}{t}, \qquad {\cal G}(t, Y)=\exp_Y t(D\exp_X)_{\exp^{-1}_XY}V.
\]
for all $Y\in B_{\delta}(X),\, t\in(0,\delta)$. Note that the above inequality is equivalent to
\[
\left|\frac{F({\cal G}(t, Y))-F(Y)}{t}-\frac{F(Y+tV)-F(Y)}{t}\right|
\leq
L_X \,\frac{d({\cal G}(t, Y), \,Y+tV)}{t}.
\]
On the other hand,  Lemma \ref{ap:0.1} implies that
\[
\lim_{t\to 0\ q\to p}\frac{d({\cal G}(t, Y),
\,Y+tV)}{t}=0,
\]
which combined with the above inequality and definitions of the generalized derivatives yields the lemma.
\end{proof}
\begin{lemma}\label{ap:1}
Let $\Omega\subset\mathbb{S}^n_{++}$ be an open convex set and
$I=\{1,...,m\}$. Let $F_{i}:M\to \mathbb{R}$ be a continuously
differentiable function on $\Omega$ for every $i\in I$ and $F:M
\to\mathbb{R}$ defined by
\[
F(X):=\max_{i\in I} F_i(X).
\]
Then $F$ is locally Lipschitz  on $\Omega$ and the following holds:
\[
\partial^{\circ} F(X)=\conv\{\grad F_i(X) : i\in I(X)\}, \qquad X\in
\Omega,
\]
where $I(X):=\{i : F_i(X)=F(X), \; i=1,...,m \}$.
\end{lemma}
\begin{proof}
Take $X\in\Omega$. From Proposition \ref{SubClarke1} it is
sufficient to prove that
\[
\partial^{\circ} F(X)\subset\conv\{\grad F_i(X) : i\in I(X)\}.
\]
Let $W\in\partial^{\circ} F(X)$. Then, definition of the set $\partial^{\circ} F(X)$ and \eqref{Rdt:0.0} provides
\[
F^{\circ}(X,V)\geq\langle W,V \rangle=tr\left(V\Psi''(X)W\right), \qquad V\in\mathbb{S}^n.
\]
where $\Psi(X)=-\ln\det X$. Combining Lemma~\ref{ap:0}
with the definition of the Euclidean generalized subdifferential of $F$ at $X$, $\partial_E^{\circ} F(X)$, we conclude
that
\[
\Psi''(X)W\in\partial_E^{\circ} F(X),
\]
On the other hand,
\[
\partial_E^{\circ} F(X)\subset\conv\{\nabla F_i(X) : i\in I(X)\},
\]
where $I(X):=\{i\in I : F_i(X)=F(X)\}$, see
\cite{Clarke1983} Proposition 2.3.12. Then, there exist a constant
$\alpha_i\geq 0$ for $i\in I(X)$ with $\sum_{i\in I(X)}\alpha_i=1$ such that
\[
\Psi''(X)W=\sum_{i\in I(X)}\alpha_i\nabla F_i(X).
\]
So, since $\Psi''(X)$ is invertible and $\grad F_i(X)=\Psi''(X)^{-1}\nabla
F_i(X)$, the above equality yields
\[
W=\sum_{i\in I(X)}\alpha_i\grad F_i(X).
\]
Hence $W\in\conv\{\grad F_i(X) : i\in I(X)\}$ and the result is
proved.
\end{proof}

\end{document}